\documentclass[a4paper,reqno,11pt]{amsart}
\usepackage[margin=25mm]{geometry}

\usepackage{amssymb}
\usepackage{eucal}
\usepackage{enumerate}
\usepackage{tikz-cd}
\usepackage{url}

\numberwithin{equation}{section}

\newtheorem{theorem}{Theorem}[section]
\newtheorem{corollary}[theorem]{Corollary}
\newtheorem{lemma}[theorem]{Lemma}
\newtheorem{proposition}[theorem]{Proposition}

\theoremstyle{definition}
\newtheorem{definition}[theorem]{Definition}
\newtheorem{remark}[theorem]{Remark}
\newtheorem{example}[theorem]{Example}

\newcommand{\C}{\mathcal{C}}
\newcommand{\D}{\mathcal{D}}
\newcommand{\E}{\mathcal{E}}
\newcommand{\F}{\mathcal{F}}
\renewcommand{\H}{\mathcal{H}}
\newcommand{\K}{\mathcal{K}}
\renewcommand{\S}{\mathcal{S}}
\newcommand{\T}{\mathcal{T}}
\newcommand{\U}{\mathcal{U}}

\newcommand{\LL}{\mathbb{L}}
\newcommand{\ZZ}{\mathbb{Z}}
\newcommand{\PP}{\mathbb{P}}

\newcommand{\e}{\epsilon}
\newcommand{\bte}{\Theta}
\newcommand{\bde}{\Delta}

\newcommand{\til}[1]{\widetilde{#1}}
\newcommand{\kk}{\Bbbk}
\newcommand{\xto}{\xrightarrow}
\newcommand{\lto}{\longrightarrow}

\newcommand{\id}{\mathrm{id}}
\newcommand{\op}{\mathrm{op}}
\newcommand{\ev}{\mathrm{ev}}
\DeclareMathOperator{\Co}{Cone}

\DeclareMathOperator{\pd}{pdim}
\DeclareMathOperator{\gldim}{gldim}
\let\mod\relax
\DeclareMathOperator{\mod}{\mathsf{mod}}
\DeclareMathOperator{\Mod}{\mathsf{Mod}}
\DeclareMathOperator{\proj}{\mathsf{proj}}
\DeclareMathOperator{\add}{\mathsf{add}}
\DeclareMathOperator{\thick}{\mathsf{thick}}
\DeclareMathOperator{\Loc}{\mathsf{Loc}}
\DeclareMathOperator{\per}{\mathsf{per}}
\DeclareMathOperator{\End}{End}
\DeclareMathOperator{\REnd}{\mathbb{R}End}
\DeclareMathOperator{\Hom}{Hom}
\DeclareMathOperator{\RHom}{\mathbb{R}Hom}
\newcommand{\dgmod}{\C_{\mathrm{dg}}}

\begin{document}
\title{Dlab--Ringel standardization in algebraic triangulated categories}
\author{Takahide Adachi}
\address{T.~Adachi: Faculty of Global and Science Studies, Yamaguchi University, 1677-1 Yoshida, Yamaguchi 753-8511, Japan}
\email{tadachi@yamaguchi-u.ac.jp}
\author{Mayu Tsukamoto}
\address{M.~Tsukamoto: Graduate School of Mathematics, Nagoya University, Furo-cho, Chikusa-ku, Nagoya, Aichi 464-8602, Japan}
\email{tsukamoto@math.nagoya-u.ac.jp}


\begin{abstract}
We introduce the notion of derived standardly stratified algebras as a connective dg analogue of standardly stratified algebras. To this end, we extend the Dlab--Ringel standardization method to algebraic triangulated categories and show that every basic proper derived standardly stratified algebra arises from this extended method. We establish Ringel duality for derived standardly stratified algebras using silting objects. More precisely, the Ringel dual is again derived standardly stratified, and the corresponding perfect derived categories are related by mutually quasi-inverse contravariant triangle equivalences. We also show that, over an algebraically closed field, a proper derived standardly stratified algebra is smooth if and only if its derived standard modules are exceptional. Finally, we prove that dg algebras given by finite acyclic graded quivers with arrows of non-positive degrees and homogeneous relations are derived standardly stratified with respect to every total order on the vertices.
\end{abstract}

\maketitle

\section{Introduction}

Quasi-hereditary algebras arose from the study of highest weight categories in the representation theory of complex semisimple Lie algebras and algebraic groups \cite{S87,CPS88}.
Standardly stratified algebras were introduced as a natural generalization of quasi-hereditary algebras \cite{CPS96,Dl96}.
Standard modules play a central role in the theory of both classes of algebras.
In the quasi-hereditary case, Ringel constructed the characteristic tilting module from the standard modules and proved that its endomorphism algebra is again quasi-hereditary \cite{R91}.
This construction gives rise to Ringel duality, which was subsequently extended to standardly stratified algebras \cite{AHLU00,Xi02}.

Dlab and Ringel introduced the notion of a standardizable set, an ordered set of objects in an abelian category satisfying properties analogous to those of standard modules \cite{DR92}.
Their standardization method constructs a quasi-hereditary algebra from such a set by using universal extensions.
Its counterpart for standardly stratified algebras was developed by Erdmann and S\'aenz \cite{ES03}.
Standardization has subsequently been studied in triangulated categories \cite{MS16}, exact categories \cite{S19} and extriangulated categories \cite{AT24}.
In \cite{AT24}, Ringel duality was also studied from the viewpoint of standardization.
In all of these works, standardization produces ordinary algebras.

The aim of this paper is to introduce and study a connective dg analogue of standardly stratified algebras from the viewpoint of the Dlab--Ringel standardization method.
This also provides a framework for studying structures analogous to those of highest weight categories in algebraic triangulated categories.

We introduce the notion of derived standardly stratified algebras.
Let $A$ be a dg algebra over a field $\kk$ and let $\D(A)$ denote the (unbounded) derived category of right dg $A$-modules. 
For a family $\S$ of objects in $\D(A)$, let $\F^{\geq 0}(\S)$ denote the extension closure of $\{S[r]\mid S\in\S,\ r\geq 0\}$.

\begin{definition}[Definition \ref{def:dss}]
Let $(\Lambda,\leq)$ be a finite totally ordered set.
Fix a connective dg algebra $A$ and a family $\PP=(P(i))_{i\in\Lambda}$ of objects in $\D(A)$ satisfying $\add\bigoplus_{i\in\Lambda}P(i)=\add A$.
The dg algebra $A$ is called a \emph{derived standardly stratified algebra} with respect to $(\PP,\leq)$ if there exist triangles in $\D(A)$
\begin{align}
K(i)\to P(i)\to\bde(i)\to K(i)[1]\quad(i\in\Lambda)\notag
\end{align}
satisfying the following conditions:
\begin{enumerate}
\item[(DS1)] $\End_{\D(A)}(\bde(i))$ is local for all $i\in\Lambda$.
\item[(DS2)] $\Hom_{\D(A)}(\bde(i),\bde(j)[r])=0$ for all $i>j$ and $r\in\ZZ$.
\item[(DS3)] $K(i)\in\F^{\geq 0}(\bde(>i))$ for all $i\in\Lambda$, where $\bde(>i):=\{\bde(j)\mid j>i\}$.
\end{enumerate}
We call $\bde(i)$ a \emph{derived standard module} of $A$ with respect to $(\PP,\leq)$.
\end{definition}

Every standardly stratified algebra, regarded as a dg algebra concentrated in degree zero, is derived standardly stratified with respect to the same order.
However, the converse does not hold even for finite-dimensional algebras concentrated in degree zero (Example \ref{ex:two_cycle}).

To construct derived standardly stratified algebras and establish Ringel duality, we extend the Dlab--Ringel standardization method to algebraic triangulated categories.
Let $\T$ be an algebraic triangulated category with a dg enhancement $\E$.
Let $\bte:=\{ \bte(i)\in \T \mid i\in\Lambda\}$ be a higher standardizable set with respect to $\leq$ (see Definition \ref{def:hss}), which is a higher analogue of a standardizable set \cite{DR92, ES03}.
Under a finiteness condition (HF), using universal extensions, we construct a basic silting object $P_{\bte}=\bigoplus_{i\in\Lambda} P_{\bte}(i)$ in the smallest thick subcategory $\thick_{\T}(\bte)$ containing $\bte$ (see Lemma \ref{lem:DR_stand} and Proposition \ref{prop:const_psilt}).

Our first main result gives a realization of derived standardly stratified algebras by the standardization method.

\begin{theorem}[Theorem \ref{thm:standardization}]
In the above setting, the following statements hold.
\begin{enumerate}
\item The dg algebra $A:=\End_{\E}(P_{\bte})$ is a proper derived standardly stratified algebra with respect to the family $(\Hom_{\E}(P_{\bte},P_{\bte}(i)))_{i\in\Lambda}$ and the order $\leq$.
\item Every basic proper derived standardly stratified algebra is quasi-isomorphic to a dg algebra obtained as in (1).
\end{enumerate}
\end{theorem}

We establish Ringel duality for derived standardly stratified algebras.
For this purpose, we construct \emph{characteristic silting objects} from derived standard modules as dg counterparts of characteristic tilting modules, and obtain the following result.

\begin{theorem}[Corollary \ref{cor:double-duality}, Theorems \ref{thm:Ringel_duality} and \ref{thm:Ringel_dss}]
Let $A$ be a basic proper derived standardly stratified algebra with respect to $(\PP,\leq)$.
Then the following statements hold.
\begin{enumerate}
\item There exist a characteristic silting object $T=\bigoplus_{i\in\Lambda}T(i)\in\per A$ and a basic proper derived standardly stratified algebra $R$ with respect to $(\PP',\leq^{\op})$ such that $R$ is quasi-isomorphic to $\REnd_{A}(T)^{\op}$.
Here $\PP'$ can be taken to be the family $(\RHom_{A}(T(i),T))_{i\in\Lambda}$ in $\per R$.
\item There exists a characteristic silting object $T'\in\per R$ such that $A$ is quasi-isomorphic to $\REnd_{R}(T')^{\op}$. Moreover, the functors
\begin{align}
\RHom_{A}(-,T):\per A\longrightarrow\per R \quad\text{and}\quad \RHom_{R}(-,T'):\per R\longrightarrow\per A \notag
\end{align}
are mutually quasi-inverse contravariant triangle equivalences.
\end{enumerate}
\end{theorem}

Over an algebraically closed field, we characterize the smoothness of derived standardly stratified algebras as follows.

\begin{theorem}[Theorem \ref{thm:dgqh}]
Assume that $\kk$ is algebraically closed.
Let $A$ be a proper derived standardly stratified algebra with respect to $(\PP,\leq)$. Denote the corresponding derived standard modules by $\bde(i)$ $(i\in\Lambda)$.
Then the following statements are equivalent.
\begin{enumerate}
\item $\bde(i)$ is exceptional for all $i\in\Lambda$.
\item $A$ is smooth.
\item $A$ has finite global dimension in the sense of \cite{H26}.
\item $\per A$ admits a full exceptional sequence.
\end{enumerate}
\end{theorem}

Motivated by this theorem, we call a derived standardly stratified algebra with respect to $(\PP,\leq)$ \emph{derived quasi-hereditary} if all the corresponding derived standard modules are exceptional.

We also give a class of derived quasi-hereditary algebras.
It is known that the path algebra of a finite acyclic quiver is quasi-hereditary with respect to every total order on the set of vertices \cite{DR89}.
Our last main result is a dg analogue of this fact that also allows relations.

\begin{theorem}[Corollary \ref{cor:acyclic}]
Let $Q$ be a finite acyclic graded quiver whose arrows have non-positive degrees.
Let $J$ be the arrow ideal of $\kk Q$, and let $I\subseteq J^{2}$ be a homogeneous ideal.
Equip $A=\kk Q/I$ with an arbitrary differential.
Then $A$ is derived quasi-hereditary with respect to every total order on the set of vertices of $Q$.
\end{theorem}

\subsection*{Notation and convention}

Throughout this paper, $\kk$ is a field, and the $\kk$-dual $\Hom_{\kk}(-,\kk)$ is denoted by $(-)^{\ast}$.
All categories and functors are $\kk$-linear, and all subcategories are full and closed under isomorphisms.
Modules and dg modules are right modules and right dg modules, respectively, unless otherwise stated.

Let $\T$ be a triangulated category with shift functor $[1]$.
For an object $M$ of $\T$, let $\add M$ denote the subcategory of $\T$ consisting of direct summands of finite direct sums of copies of $M$.
For a class $\S$ of objects in $\T$, let $\thick_{\T}(\S)$ denote the smallest thick subcategory of $\T$ containing $\S$.
Let $\F(\S)$ denote the smallest extension-closed subcategory of $\T$ containing $\S$ and the zero object, and put
\begin{align}
\F^{\geq 0}(\S):=\F(\{S[r]\mid S\in\S,\ r\geq 0\}) \quad\text{and}\quad \F^{\leq 0}(\S):=\F(\{S[r]\mid S\in\S,\ r\leq 0\}).\notag
\end{align}
For an object $M$ of $\T$, we write $\thick_{\T}(M):=\thick_{\T}(\{M\})$, $\F(M):=\F(\{M\})$, $\F^{\geq 0}(M):=\F^{\geq 0}(\{M\})$ and $\F^{\leq 0}(M):=\F^{\leq 0}(\{M\})$.

All gradings are cohomological.
All opposite dg algebras and dg categories are taken with the usual Koszul sign convention.
For a dg algebra $A$, let $\dgmod(A)$ denote the dg category of right dg $A$-modules, whose Hom complexes consist of graded $A$-linear maps with the usual differential.
Its homotopy category is denoted by $\K(A):=\H^{0}(\dgmod(A))$, and the derived category of $A$ is denoted by $\D(A)$.
Let $\per A:=\thick_{\D(A)}(A)$.
We say that $A$ is \emph{connective} if $H^{r}(A)=0$ for all $r>0$, \emph{proper} if $\bigoplus_{r\in\ZZ}H^{r}(A)$ is finite-dimensional, and \emph{smooth} if $A\in\per(A^{\op}\otimes_{\kk}A)$, where $A$ is regarded as a dg $A$-bimodule.
Two dg algebras $A$ and $B$ are \emph{derived Morita equivalent} if there exists a dg $A$-$B$-bimodule $X$ such that $-\otimes^{\LL}_{A}X\colon\D(A)\to\D(B)$ is a triangle equivalence.

For a quiver, we compose arrows from left to right.
For example, the composition of arrows $a\colon1\to2$ and $b\colon2\to3$ is written as $ab$.

\section{Dlab--Ringel standardization method}\label{section:DR_std_meth}

In this section, we formulate the Dlab--Ringel standardization method in the setting of algebraic triangulated categories. We construct silting objects from higher standardizable sets and establish a dg version of Ringel duality.

\subsection{Universal extension}\label{subsection:univ_ext}

Throughout this subsection, let $\T$ be a $\kk$-linear triangulated category.
We construct presilting objects from higher standardizable sets using universal extensions.
We also give a construction using universal coextensions and introduce characteristic silting objects.

We begin by introducing the notion of a higher standardizable set in $\T$.

\begin{definition}\label{def:hss}
Let $(\Lambda,\leq)$ be a finite totally ordered set.
A set $\bte=\{\bte(i) \in \T \mid i\in\Lambda\}$ is called a \emph{higher standardizable set} with respect to $\leq$ if the following conditions hold.
\begin{enumerate}
\item[(S1)] $\End_{\T}(\bte(i))$ is local for all $i\in\Lambda$.
\item[(S2)] $\Hom_{\T}(\bte(i),\bte(i)[r])=0$ for all $i\in\Lambda$ and $r>0$.
\item[(S3)] $\Hom_{\T}(\bte(i),\bte(j)[r])=0$ for all $i>j$ and $r\in\ZZ$.
\end{enumerate}
If moreover $\thick_{\T}(\bte)=\T$, then $\bte$ is said to be \emph{full}. 
We say that $\bte$ satisfies (HF) if 
\begin{align}
\dim_{\kk}\bigoplus_{r\in\ZZ} \Hom_{\T}(\bte(i),\bte(j)[r])<\infty\notag
\end{align}
for all $i,j\in \Lambda$.
\end{definition}

We give a typical example of higher standardizable sets.

\begin{example}\label{ex:exceptional}
An object $E$ of $\T$ is said to be \emph{exceptional} if $\End_{\T}(E)$ is a division algebra and $\Hom_{\T}(E,E[r])=0$ for all $r\neq 0$.
A sequence $(E_{1},\ldots,E_{n})$ of exceptional objects of $\T$ is called an \emph{exceptional sequence} if $\Hom_{\T}(E_{i},E_{j}[r])=0$ for all $i>j$ and $r\in\ZZ$.
It is said to be \emph{full} if $\thick_{\T}(\{E_{1},\ldots,E_{n}\})=\T$.
Every exceptional sequence $(E_{1},\ldots,E_{n})$ is a higher standardizable set with respect to the natural order on $\{1,\ldots,n\}$.
\end{example}

Fix a higher standardizable set $\bte=\{\bte(i)\in \T \mid i\in\Lambda\}$ with respect to $\leq$.
For each $i\in\Lambda$, put $\bte(>i):=\{\bte(j)\mid j>i\}$ and $\bte(<i):=\{\bte(j)\mid j<i\}$.
Clearly, $\F^{\geq 0}(\bte(>i))$ and $\F^{\geq 0}(\bte)$ are closed under positive shifts, while $\F^{\leq 0}(\bte(<i))$ and $\F^{\leq 0}(\bte)$ are closed under negative shifts.
By repeated applications of \cite[Lemma 2.6(a)]{IY18}, these subcategories are closed under direct summands.
Note that $\thick_{\T}(\bte)$ is idempotent complete.
Indeed, let $\T^{\omega}$ be an idempotent completion of $\T$. By \cite{BS01}, $\T^{\omega}$ is a triangulated category, and $\T$ can be viewed as an extension-closed full triangulated subcategory of $\T^{\omega}$. Then we can check that $\thick_{\T}(\bte)=\thick_{\T^{\omega}}(\bte)$ by (S1)--(S3). Since $\thick_{\T^{\omega}}(\bte)$ is idempotent complete, so is $\thick_{\T}(\bte)$.

Throughout this section, we always assume that $\bte$ satisfies (HF).
By a standard d\'evissage argument, $\bigoplus_{r\in\ZZ}\Hom_{\T}(X,Y[r])$ is finite-dimensional for all $X,Y\in\thick_{\T}(\bte)$.
In particular, $\thick_{\T}(\bte)$ is Hom-finite, and hence a Krull--Schmidt category.

Using universal extensions, we construct a presilting object in $\T$ from a higher standardizable set.
This construction is an analogue of that in \cite[Appendix A]{CC26}.
Recall the notion of universal extension.
Let $X,Y\in\T$ and $r>0$.
Assume that $V:=\Hom_{\T}(X,Y[r])$ is finite-dimensional over $\kk$.
Consider a triangle
\begin{align}\label{univ_tri}
V^{\ast}\otimes_{\kk}Y[r-1]\to E\to X \xto{\e}V^{\ast}\otimes_{\kk}Y[r],
\end{align}
where $\e$ corresponds to $\id_{V}$ under the canonical isomorphisms
\begin{align}
\Hom_{\T}(X,V^{\ast}\otimes_{\kk}Y[r]) \cong V^{\ast}\otimes_{\kk}V \cong\End_{\kk}(V).\notag
\end{align}
We call this triangle a \emph{universal extension}.

\begin{lemma}\label{lem:high_univ_ext}
The connecting map
\begin{align}
\delta:\Hom_{\T}(V^{\ast}\otimes_{\kk}Y[r],Y[r])
\to\Hom_{\T}(X,Y[r])\notag
\end{align}
is surjective. In particular, if $\Hom_{\T}(Y,Y[1])=0$, then $\Hom_{\T}(E,Y[r])=0$.
\end{lemma}
\begin{proof}
Under the canonical isomorphism $\Hom_{\T}(V^{\ast}\otimes_{\kk}Y[r],Y[r])\cong V\otimes_{\kk}\End_{\T}(Y)$, the map $\delta$ is induced by composition.
It admits a section $V\to V\otimes_{\kk}\End_{\T}(Y)$ ($v\mapsto v\otimes\id_{Y}$).
Thus $\delta$ is surjective.
By applying $\Hom_{\T}(-,Y)$ to the triangle \eqref{univ_tri}, the last assertion follows from $\Hom_{\T}(Y,Y[1])=0$.
\end{proof}

We construct the following triangles associated with $\bte$ using universal extensions.

\begin{lemma}\label{lem:DR_stand}
For each $i\in\Lambda$, there exists a triangle
\begin{align}
K(i)\to P(i)\to \bte(i)\to K(i)[1]\notag
\end{align}
such that $K(i)\in\F^{\geq 0}(\bte(>i))$ and $\Hom_{\T}(P(i),Y[r])=0$ for all $Y\in\F^{\geq 0}(\bte)$ and $r>0$.
\end{lemma}
\begin{proof}
Write $n=|\Lambda|$ and identify $\Lambda$ with $[1,n]$ via the order-preserving bijection.
Fix $i\in[1,n]$.
For each $m\in[i,n]$, we inductively construct an object $P(i,m)$ together with a triangle
\begin{align}
K(i,m)\to P(i,m)\to\bte(i)\to K(i,m)[1]\notag
\end{align}
such that $K(i,m)\in\F^{\geq 0}(\bte(i,m])$ and $\Hom_{\T}(P(i,m),\bte(j)[r])=0$ for all $j\in[1,m]$ and $r>0$, where
\begin{align}
\F^{\geq 0}(\bte(i,m]):= \F(\{\bte(j)[r]\mid j\in (i,m],r\geq 0\}).\notag
\end{align}

If $m=i$, we take $P(i,i):=\bte(i)$ and $K(i,i):=0$. The assertion follows from (S2) and (S3).
Assume that $m>i$. By the induction hypothesis, $P(i,m-1)$ and $K(i,m-1)$ have been constructed.
Let
\begin{align}
N:=\max\bigl( \{0\}\cup \{r>0\mid\Hom_{\T}(P(i,m-1),\bte(m)[r])\neq 0\} \bigr).\notag
\end{align}
By the condition (HF), the integer $N$ is well-defined.

For $t\in [0,N]$, we inductively construct triangles
\begin{align}
K_{t}\to X_{t}\to \bte(i)\to K_{t}[1]\notag
\end{align}
satisfying $K_{t}\in\F^{\geq 0}(\bte(i,m])$, $\Hom_{\T}(X_{t},\bte(m)[r])=0$ for all $r>t$ and $\Hom_{\T}(X_{t},\bte(j)[r])=0$ for all $j<m$ and $r>0$.
Let $X_{N}:=P(i,m-1)$ and $K_{N}:=K(i,m-1)$.
Assume that $t<N$ and the triangle for $t+1$ has been constructed.
Put $V_{t}:=\Hom_{\T}(X_{t+1},\bte(m)[t+1])$.
By $X_{t+1}\in \thick_{\T}(\bte)$, the $\kk$-vector space $V_{t}$ is finite-dimensional.
Consider the universal extension
\begin{align}
V_{t}^{\ast}\otimes_{\kk}\bte(m)[t] \to X_{t}\to X_{t+1}\xto{\e_{t}}V_{t}^{\ast}\otimes_{\kk}\bte(m)[t+1].\notag
\end{align}
By the octahedral axiom, we have triangles
\begin{align}
K_{t}\to X_{t}\to \bte(i)\to K_{t}[1]\notag
\end{align}
and
\begin{align}
V_{t}^{\ast}\otimes_{\kk}\bte(m)[t]\to K_{t}\to K_{t+1}\to V_{t}^{\ast}\otimes_{\kk}\bte(m)[t+1].\notag
\end{align}
Since $\F^{\geq 0}(\bte(i,m])$ is closed under extensions, we have $K_{t}\in \F^{\geq 0}(\bte(i,m])$.
For $j<m$, the induction hypothesis and condition (S3) give $\Hom_{\T}(X_{t+1},\bte(j)[r])=0$ for all $r>0$ and $\Hom_{\T}(V_{t}^{\ast}\otimes_{\kk}\bte(m)[t],\bte(j)[r])=0$ for all $r\in\ZZ$, respectively.
Hence $\Hom_{\T}(X_{t},\bte(j)[r])=0$ for all $j<m$ and $r>0$.
For $j=m$ and $r>t+1$, the induction hypothesis and (S2) give $\Hom_{\T}(X_{t+1},\bte(m)[r])=0$ and $\Hom_{\T}(V_{t}^{\ast}\otimes_{\kk}\bte(m)[t], \bte(m)[r])=0$.
Consequently, $\Hom_{\T}(X_{t},\bte(m)[r])=0$ for all $r>t+1$.
For $j=m$ and $r=t+1$, Lemma \ref{lem:high_univ_ext} and (S2) give $\Hom_{\T}(X_{t},\bte(m)[t+1])=0$.

Let $P(i,m):=X_{0}$ and $K(i,m):=K_{0}$. Then the required conditions hold, completing the induction on $m$.
Finally, let $P(i):=P(i,n)$ and $K(i):=K(i,n)$.
By the definition of $\F^{\geq 0}(\bte)$, we obtain $\Hom_{\T}(P(i),Y[r])=0$ for all $Y\in\F^{\geq 0}(\bte)$ and $r>0$.
\end{proof}

We establish the silting property of the object constructed above.
Recall that an object $M$ of $\T$ is called \emph{presilting} if $\Hom_{\T}(M,M[r])=0$ for all $r>0$, and \emph{silting} if moreover $\thick_{\T}(M)=\T$ (see \cite{KV88,AI12}).

\begin{proposition}\label{prop:const_psilt}
For each $i\in\Lambda$, consider any triangle
\begin{align}\label{eq:psilt_tri}
K_{i}\to P_{i}\xto{\alpha_{i}}\bte(i)\to K_{i}[1]
\end{align}
such that $K_{i}\in\F^{\geq 0}(\bte(>i))$ and $\Hom_{\T}(P_{i},Y[r])=0$ for all $Y\in\F^{\geq 0}(\bte)$ and $r>0$.
Let $P:=\bigoplus_{i\in\Lambda}P_{i}$.
Then the following statements hold.
\begin{enumerate}
\item $P$ is a presilting object in $\T$.
\item $\thick_{\T}(P)=\thick_{\T}(\bte)$. In particular, $P$ is a silting object in $\thick_{\T}(\bte)$.
\item If $\bte$ is full, then $P$ is a silting object in $\T$.
\end{enumerate}
\end{proposition}
\begin{proof}
(1) Since $K_{i}\in\F^{\geq 0}(\bte(>i))$, the triangle \eqref{eq:psilt_tri} gives $P_{i}\in\F^{\geq 0}(\bte)$ for all $i\in\Lambda$.
Thus $P$ is presilting.

(2) By \eqref{eq:psilt_tri}, we have $P_{i}\in\thick_{\T}(\bte)$ for all $i\in\Lambda$.
For the maximal element $m\in\Lambda$, we have $K_{m}=0$, and hence $\bte(m)\cong P_{m}$.
Since $K_{i}\in\F^{\geq 0}(\bte(>i))$, descending induction on $i$ yields $\bte(i)\in\thick_{\T}(P)$ for all $i\in\Lambda$.
Thus $\thick_{\T}(P)=\thick_{\T}(\bte)$.
By (1), $P$ is silting in $\thick_{\T}(\bte)$.

(3) Since $\bte$ is full, we have $\T=\thick_{\T}(\bte)$. 
Thus the assertion follows from (2).
\end{proof}

Using universal coextensions, we give the counterparts of Lemma \ref{lem:DR_stand} and Proposition \ref{prop:const_psilt}.

\begin{lemma}\label{lem:high_univ_coext}
For each $i\in\Lambda$, there exists a triangle
\begin{align}
\bte(i)\to I(i)\to C(i)\to\bte(i)[1]\notag
\end{align}
such that $C(i)\in\F^{\leq 0}(\bte(<i))$ and $\Hom_{\T}(X,I(i)[r])=0$ for all $X\in\F^{\leq 0}(\bte)$ and $r>0$.
\end{lemma}

\begin{proposition}\label{prop:coext_silting}
For each $i\in\Lambda$, consider any triangle
\begin{align}\label{eq:coext_tri}
\bte(i)\xto{\beta_{i}}I_{i}\to C_{i}\to\bte(i)[1]
\end{align}
such that $C_{i}\in\F^{\leq 0}(\bte(<i))$ and $\Hom_{\T}(X,I_{i}[r])=0$ for all $X\in\F^{\leq 0}(\bte)$ and $r>0$.
Let $I:=\bigoplus_{i\in\Lambda}I_{i}$.
Then the following statements hold.
\begin{enumerate}
\item $I$ is a presilting object in $\T$.
\item $\thick_{\T}(I)=\thick_{\T}(\bte)$. In particular, $I$ is a silting object in $\thick_{\T}(\bte)$.
\item If $\bte$ is full, then $I$ is a silting object in $\T$.
\end{enumerate}
\end{proposition}

\begin{remark}\label{rem:DR-minimal}
Recall that $\thick_{\T}(\bte)$ is a Krull--Schmidt category.
Consider any triangles \eqref{eq:psilt_tri} as in Proposition \ref{prop:const_psilt}.
Since $\F^{\geq 0}(\bte(>i))$ is closed under direct summands, we may assume that each $\alpha_{i}$ is right minimal by taking its right minimal version.
Then $P_{i}$ is indecomposable by (S1) and (S3).
Moreover, any two such triangles with right minimal morphisms to $\bte(i)$ are isomorphic.
Thus $P_{i}$ is uniquely determined up to isomorphism by $\bte$ and the chosen order.
In particular, the object $P_{\bte}:=\bigoplus_{i\in\Lambda}P_{\bte}(i)$, where $P_{\bte}(i):=P_{i}$, is basic and uniquely determined up to isomorphism.
Similarly, for any triangles \eqref{eq:coext_tri} as in Proposition \ref{prop:coext_silting}, we may assume that each $\beta_{i}$ is left minimal, and then each $I_{i}$ is indecomposable and uniquely determined up to isomorphism.
In particular, the object $I_{\bte}:=\bigoplus_{i\in\Lambda}I_{\bte}(i)$, where $I_{\bte}(i):=I_{i}$, is basic and uniquely determined up to isomorphism.
\end{remark}

Characteristic tilting modules play a central role in Ringel duality for quasi-hereditary algebras.
Motivated by the preceding construction, we introduce the following silting analogue.

\begin{definition}
An object $T$ in $\thick_{\T}(\bte)$ is called a \emph{characteristic silting object} associated with $\bte$ if $T\cong\bigoplus_{i\in\Lambda}I_{i}$ for some triangles \eqref{eq:coext_tri} as in Proposition \ref{prop:coext_silting}.
\end{definition}

By Lemma \ref{lem:high_univ_coext}, the object $\bigoplus_{i\in\Lambda}I(i)$ is a characteristic silting object associated with $\bte$.

\subsection{Ringel duality in algebraic triangulated categories}\label{subsection:Ringel_alg}

In this subsection, we investigate the relationship between higher standardizable sets in algebraic triangulated categories and the dg algebras arising from their associated silting objects.
We establish a dg version of Ringel duality between these dg algebras.

Throughout this subsection, let $\E$ be a $\kk$-linear pretriangulated dg category.
Recall that the homotopy category $\T:=\H^{0}(\E)$ has the same objects as $\E$ and $\Hom_{\T}(X,Y)=H^{0}\Hom_{\E}(X,Y)$.
Since $\E$ is pretriangulated, $\T$ is a triangulated category with shift functor $[1]$, and we have
\begin{align}
H^{r}\Hom_{\E}(X,Y)\cong\Hom_{\T}(X,Y[r])\notag
\end{align}
for all $X,Y\in\T$ and $r\in\ZZ$.
For a dg functor $\Phi:\C\to\C'$, the induced functor
\begin{align}
\H^{0}(\Phi):\H^{0}(\C)\longrightarrow\H^{0}(\C')\notag
\end{align}
is defined by $\H^{0}(\Phi)(X)=\Phi(X)$ on objects and $\H^{0}(\Phi)([f])=[\Phi(f)]$ on morphisms, where $f\in Z^{0}\Hom_{\C}(X,Y)$.
If $\C$ and $\C'$ are pretriangulated, then $\H^{0}(\Phi)$ is a triangle functor.

Let $\bte=\{\bte(i)\in\T\mid i\in\Lambda\}$ be a higher standardizable set with respect to $\leq$ satisfying (HF).
By Lemma \ref{lem:DR_stand} and Remark \ref{rem:DR-minimal}, we can choose triangles
\begin{align}\label{eq:P_tri}
K_{\bte}(i)\to P_{\bte}(i)\to \bte(i)\to K_{\bte}(i)[1]\quad(i\in\Lambda)
\end{align}
as in Proposition \ref{prop:const_psilt} such that $P_{\bte}=\bigoplus_{i\in\Lambda}P_{\bte}(i)$ is a basic silting object in $\thick_{\T}(\bte)$, where each $P_{\bte}(i)$ is indecomposable.
By Lemma \ref{lem:high_univ_coext} and Remark \ref{rem:DR-minimal}, we can choose triangles
\begin{align}\label{eq:I_tri}
\bte(i)\to I_{\bte}(i)\to C_{\bte}(i)\to\bte(i)[1]\quad(i\in\Lambda)
\end{align}
as in Proposition \ref{prop:coext_silting} such that $I_{\bte}=\bigoplus_{i\in\Lambda}I_{\bte}(i)$ is a basic silting object in $\thick_{\T}(\bte)$, where each $I_{\bte}(i)$ is indecomposable.

Let $A:=\End_{\E}(P_{\bte})$.
Since $P_{\bte}$ is presilting and $\bte$ satisfies (HF), the dg algebra $A$ is connective and proper.
Hence $\per A$ is a Krull--Schmidt category.
For any $Y\in\T$, the complex $\Hom_{\E}(P_{\bte},Y)$ is naturally a right dg $A$-module via composition.
Since $\E$ and $\dgmod(A)$ are pretriangulated, the dg functor
\begin{align}
\Hom_{\E}(P_{\bte},-):\E\longrightarrow\dgmod(A)\notag
\end{align}
induces a triangle functor
\begin{align}
\T=\H^{0}(\E) \xrightarrow{\H^{0}(\Hom_{\E}(P_{\bte},-))} \H^{0}(\dgmod(A))=\K(A)\xrightarrow{\pi}\D(A),\notag
\end{align}
where $\pi$ is the canonical localization.
We denote this composite by $F$.
Thus, for $Y\in\T$, the object $FY$ of $\D(A)$ is represented by the dg $A$-module $\Hom_{\E}(P_{\bte},Y)$.

Let $R:=\End_{\E}(I_{\bte})^{\op}$.
Since $I_{\bte}$ is presilting and $\bte$ satisfies (HF), the dg algebra $R$ is connective and proper, and hence $\per R$ is a Krull--Schmidt category.
For any $Y\in\T$, the complex $\Hom_{\E}(Y,I_{\bte})$ is naturally a right dg $R$-module via $f\cdot r:=(-1)^{|f||r|}r\circ f$ for homogeneous $f\in\Hom_{\E}(Y,I_{\bte})$ and $r\in R$.
The contravariant dg functor
\begin{align}
\Hom_{\E}(-,I_{\bte}):\E\longrightarrow\dgmod(R)\notag
\end{align}
induces a contravariant triangle functor
\begin{align}
\T=\H^{0}(\E) \xrightarrow{\H^{0}(\Hom_{\E}(-,I_{\bte}))} \H^{0}(\dgmod(R))=\K(R)\xrightarrow{\pi}\D(R),\notag
\end{align}
where $\pi$ is the canonical localization.
We denote this composite by $G$.
Thus, for $Y\in\T$, the object $GY$ of $\D(R)$ is represented by the dg $R$-module $\Hom_{\E}(Y,I_{\bte})$.

The following two lemmas play an important role in this subsection.
A dg $A$-module $U$ is called \emph{homotopically projective} (or \emph{$K$-projective}) if $\Hom_{\K(A)}(U,V)=0$ for every acyclic dg $A$-module $V$.

\begin{lemma}\label{lem:Yoneda}
The following statements hold.
\begin{enumerate}
\item For every $X\in\thick_{\T}(P_{\bte})$, the dg $A$-module $\Hom_{\E}(P_{\bte},X)$ is homotopically projective.
\item $F(\thick_{\T}(P_{\bte}))\subseteq\per A$.
\item For all $X\in\thick_{\T}(P_{\bte})$ and $Y\in\T$, the morphism of complexes
\begin{align}
\Hom_{\E}(X,Y)\lto\Hom_{\dgmod(A)}(\Hom_{\E}(P_{\bte},X),\Hom_{\E}(P_{\bte},Y))\notag
\end{align}
induced by $\Hom_{\E}(P_{\bte},-)$ is a quasi-isomorphism.
\item $F:\thick_{\T}(P_{\bte})\to\per A$ is a fully faithful triangle functor.
\end{enumerate}
\end{lemma}

This is well known (see \cite{K94, K06}); we include a proof for the convenience of the reader.

\begin{proof}
Let $\til{F}:=\H^{0}(\Hom_{\E}(P_{\bte},-)):\T\to\K(A)$.
Since $\til{F}$ is a triangle functor and $\til{F}(P_{\bte})=A$, we have
\begin{align}\label{eq:tilFthick}
\til{F}(\thick_{\T}(P_{\bte}))\subseteq\thick_{\K(A)}(A).
\end{align}

(1) Since the homotopically projective dg $A$-modules form a thick subcategory of $\K(A)$ containing $A$, the assertion follows from \eqref{eq:tilFthick}.

(2) Since $\pi$ is a triangle functor, we have
\begin{align}
\pi(\thick_{\K(A)}(A))\subseteq\thick_{\D(A)}(A)=\per A.\notag
\end{align}
Combining this with the inclusion \eqref{eq:tilFthick}, we obtain $F(\thick_{\T}(P_{\bte}))\subseteq\per A$.

(3) For $X,Y\in\T$, denote by $\phi_{X,Y}$ the morphism of complexes in (3).
By definition, $H^{0}(\phi_{X,Y})$ is the map $\Hom_{\T}(X,Y)\to\Hom_{\K(A)}(\til{F}X,\til{F}Y)$ induced by $\til{F}$.
For every $Y\in\T$, the dg Yoneda lemma \cite[1.2]{K94} gives an isomorphism of complexes
\begin{align}
\phi_{P_{\bte},Y}:\Hom_{\E}(P_{\bte},Y)\lto\Hom_{\dgmod(A)}(A,\Hom_{\E}(P_{\bte},Y))\quad(u\mapsto u\circ-).\notag
\end{align}
In particular, $H^{0}(\phi_{P_{\bte},Y})$ is bijective.
Since $\til{F}$ is a triangle functor, $H^{0}(\phi_{X,Y})$ is bijective for all $X\in\thick_{\T}(P_{\bte})$ and $Y\in\T$ by a standard d\'evissage argument in the first variable.

Let $r\in\ZZ$.
Since $\E$ is pretriangulated, there is a homotopy equivalence
\begin{align}
\theta_{Z}:\Hom_{\E}(Z,Y[r])\lto\Hom_{\E}(Z,Y)[r]\notag
\end{align}
natural in $Z\in\E$.
In particular, $\theta_{P_{\bte}}:\til{F}(Y[r])\to(\til{F}Y)[r]$ is a homotopy equivalence of right dg $A$-modules, and hence the map
\begin{align}
\theta_{P_{\bte}}\circ-:\Hom_{\dgmod(A)}(\til{F}X,\til{F}(Y[r]))\lto\Hom_{\dgmod(A)}(\til{F}X,\til{F}Y)[r]\notag
\end{align}
is a homotopy equivalence, where we identify $\Hom_{\dgmod(A)}(\til{F}X,(\til{F}Y)[r])$ with $\Hom_{\dgmod(A)}(\til{F}X,\til{F}Y)[r]$.
By the naturality of $\theta_{Z}$, we have $(\theta_{P_{\bte}}\circ-)\circ\phi_{X,Y[r]}=\phi_{X,Y}[r]\circ\theta_{X}$.
Applying $H^{0}$, we obtain
\begin{align}
H^{r}(\phi_{X,Y})\circ H^{0}(\theta_{X})=H^{0}(\phi_{X,Y}[r])\circ H^{0}(\theta_{X})=H^{0}(\theta_{P_{\bte}}\circ-)\circ H^{0}(\phi_{X,Y[r]}).\notag
\end{align}
Since $H^{0}(\theta_{X})$, $H^{0}(\theta_{P_{\bte}}\circ-)$ and $H^{0}(\phi_{X,Y[r]})$ are bijective, so is $H^{r}(\phi_{X,Y})$.
Therefore $\phi_{X,Y}$ is a quasi-isomorphism for all $X\in\thick_{\T}(P_{\bte})$ and $Y\in\T$.

(4) By (2), $F$ restricts to a triangle functor $\thick_{\T}(P_{\bte})\to\per A$.
Let $X,X'\in\thick_{\T}(P_{\bte})$.
By (3), $\til{F}$ induces a bijection $\Hom_{\T}(X,X')\to\Hom_{\K(A)}(\til{F}X,\til{F}X')$.
By (1) and \cite[3.1 and 4.1]{K94}, the localization $\pi$ induces a bijection $\Hom_{\K(A)}(\til{F}X,\til{F}X')\to\Hom_{\D(A)}(FX,FX')$.
Thus $F$ is fully faithful on $\thick_{\T}(P_{\bte})$.
\end{proof}

\begin{lemma}\label{lem:Yoneda_op}
The following statements hold.
\begin{enumerate}
\item For every $X\in\thick_{\T}(I_{\bte})$, the dg $R$-module $\Hom_{\E}(X,I_{\bte})$ is homotopically projective.
\item $G(\thick_{\T}(I_{\bte}))\subseteq\per R$.
\item For all $X\in\thick_{\T}(I_{\bte})$ and $Y\in\T$, the morphism of complexes
\begin{align}
\Hom_{\E}(Y,X)\lto\Hom_{\dgmod(R)}(\Hom_{\E}(X,I_{\bte}),\Hom_{\E}(Y,I_{\bte}))\notag
\end{align}
induced by $\Hom_{\E}(-,I_{\bte})$ is a quasi-isomorphism.
\item $G:\thick_{\T}(I_{\bte})\to\per R$ is a fully faithful contravariant triangle functor.
\end{enumerate}
\end{lemma}
\begin{proof}
The assertions follow by the same argument as in the proof of Lemma \ref{lem:Yoneda}.
\end{proof}

For each $i\in\Lambda$, let $\bde(i):=F\bte(i)$ and $\bde:=\{\bde(i)\mid i\in\Lambda\}$.
Let $T(i):=FI_{\bte}(i)$ and $T:=FI_{\bte}=\bigoplus_{i\in\Lambda}T(i)$.
We obtain the following result from Lemma \ref{lem:Yoneda}.

\begin{proposition}\label{prop:fullyfaithful}
The following statements hold.
\begin{enumerate}
\item The functor $F$ induces a triangle equivalence
\begin{align}
F:\thick_{\T}(\bte)\longrightarrow\per A.\notag
\end{align}
\item The functor $F$ induces equivalences of categories
\begin{align}
\F^{\geq 0}(\bte)\lto\F^{\geq 0}(\bde) \quad\text{and}\quad \F^{\leq 0}(\bte)\lto\F^{\leq 0}(\bde).\notag
\end{align}
\item The set $\bde$ is a higher standardizable set in $\per A$ with respect to $\leq$ and satisfies (HF).
\end{enumerate}
\end{proposition}
\begin{proof}
(1) By Lemma \ref{lem:Yoneda}(4), the essential image of $F:\thick_{\T}(P_{\bte})\to\per A$ is a triangulated subcategory of $\per A$ containing $A=FP_{\bte}$.
Since $F$ is fully faithful and $\thick_{\T}(P_{\bte})=\thick_{\T}(\bte)$ is idempotent complete, the essential image is closed under direct summands.
Hence it is a thick subcategory containing $A$, and thus equals $\per A$.

(2) By (1), $F$ is a triangle equivalence with $F(\bte(i)[r])\cong\bde(i)[r]$ for all $i\in\Lambda$ and $r\in\ZZ$.
Since a triangle equivalence sends the extension closure of a class of objects to the extension closure of its image, we have the assertion.

(3) By (1), we have an isomorphism
\begin{align}
\Hom_{\D(A)}(\bde(i),\bde(j)[r]) \cong\Hom_{\T}(\bte(i),\bte(j)[r])\notag
\end{align}
for all $i,j\in\Lambda$ and $r\in\ZZ$.
Since $\bte$ satisfies (S1), (S2), (S3) and (HF), so does $\bde$.
\end{proof}

By Proposition \ref{prop:fullyfaithful}(1) and (3), the set $\bde$ is a higher standardizable set in $\per A$ satisfying (HF), and $\thick_{\D(A)}(\bde)=\per A$.
Hence, applying the results of Subsection \ref{subsection:univ_ext} to $\bde$ in place of $\bte$, we obtain the basic silting objects $P_{\bde}=\bigoplus_{i\in\Lambda}P_{\bde}(i)$ and $I_{\bde}=\bigoplus_{i\in\Lambda}I_{\bde}(i)$ in $\per A$, where $P_{\bde}(i)$ and $I_{\bde}(i)$ are indecomposable.

\begin{proposition}\label{prop:characteristic-silting}
For each $i\in\Lambda$, we have $FP_{\bte}(i)\cong P_{\bde}(i)$ and $FI_{\bte}(i)\cong I_{\bde}(i)$.
In particular, we have $A\cong P_{\bde}$ and $T\cong I_{\bde}$, and hence $T$ is a basic characteristic silting object associated with $\bde$.
\end{proposition}
\begin{proof}
Let $i\in\Lambda$.
Applying $F$ to the triangle \eqref{eq:I_tri}, we obtain a triangle
\begin{align}
\bde(i)\to T(i)\to FC_{\bte}(i)\to\bde(i)[1]\notag
\end{align}
with $FC_{\bte}(i)\in\F^{\leq 0}(\bde(<i))$, since $F$ is a triangle functor.
Moreover, for any $X\in\F^{\leq 0}(\bde)$, Proposition \ref{prop:fullyfaithful}(2) gives $Y\in\F^{\leq 0}(\bte)$ such that $FY\cong X$.
Hence, by Proposition \ref{prop:fullyfaithful}(1), we have
\begin{align}
\Hom_{\D(A)}(X,T(i)[r]) \cong\Hom_{\T}(Y,I_{\bte}(i)[r])=0\notag
\end{align}
for all $r>0$, where the last equality holds since \eqref{eq:I_tri} is as in Proposition \ref{prop:coext_silting}.
Thus these triangles are as in Proposition \ref{prop:coext_silting} for $\bde$.
Since each $T(i)$ is indecomposable by Proposition \ref{prop:fullyfaithful}(1), Remark \ref{rem:DR-minimal} gives $T(i)\cong I_{\bde}(i)$.
Similarly, applying $F$ to \eqref{eq:P_tri}, we obtain triangles as in Proposition \ref{prop:const_psilt} for $\bde$, and hence $FP_{\bte}(i)\cong P_{\bde}(i)$.
\end{proof}

For each $i\in\Lambda$, let $\bde'(i):=G\bte(i)$ and $\bde':=\{\bde'(i)\mid i\in\Lambda\}$.
Let $T':=GP_{\bte}$.
We give the contravariant counterparts of Propositions \ref{prop:fullyfaithful} and \ref{prop:characteristic-silting}.

\begin{proposition}\label{prop:dualfullyfaithful}
The following statements hold.
\begin{enumerate}
\item The functor $G$ induces a contravariant triangle equivalence
\begin{align}
G:\thick_{\T}(\bte)\longrightarrow\per R.\notag
\end{align}
\item The functor $G$ induces contravariant equivalences of categories
\begin{align}
\F^{\geq 0}(\bte)\lto\F^{\leq 0}(\bde') \quad\text{and}\quad \F^{\leq 0}(\bte)\lto\F^{\geq 0}(\bde').\notag
\end{align}
\item The set $\bde'$ is a higher standardizable set in $\per R$ with respect to the opposite order $\leq^{\op}$ and satisfies (HF).
\end{enumerate}
\end{proposition}

By Proposition \ref{prop:dualfullyfaithful}(1) and (3), the set $\bde'$ is a higher standardizable set in $\per R$ with respect to $\leq^{\op}$ satisfying (HF), and $\thick_{\D(R)}(\bde')=\per R$.
Hence we obtain the basic silting objects $P_{\bde'}=\bigoplus_{i\in\Lambda}P_{\bde'}(i)$ and $I_{\bde'}=\bigoplus_{i\in\Lambda}I_{\bde'}(i)$ in $\per R$ with respect to $\leq^{\op}$ in the same way, where $P_{\bde'}(i)$ and $I_{\bde'}(i)$ are indecomposable.

\begin{proposition}\label{prop:dual-silting}
For each $i\in\Lambda$, we have $GI_{\bte}(i)\cong P_{\bde'}(i)$ and $GP_{\bte}(i)\cong I_{\bde'}(i)$.
In particular, we have $R\cong P_{\bde'}$ and $T'\cong I_{\bde'}$, and hence $T'$ is a basic characteristic silting object associated with $\bde'$.
\end{proposition}

To establish Ringel duality, we compute $\REnd_{A}(T)$ and $\REnd_{R}(T')$.
The complex $M:=\Hom_{\E}(P_{\bte},I_{\bte})$ is a right dg $A$-module representing $T$.
The complex $M$ also represents $T'$ as a right dg $R$-module via composition.
Since $I_{\bte}\in\thick_{\T}(P_{\bte})$ and $P_{\bte}\in\thick_{\T}(I_{\bte})$, the dg module $M$ is homotopically projective both over $A$ and over $R$ by Lemmas \ref{lem:Yoneda}(1) and \ref{lem:Yoneda_op}(1).
Hence $\REnd_{A}(T)$ and $\REnd_{R}(T')$ can be computed by the ordinary dg endomorphism algebras of $M$ over $A$ and over $R$, respectively.

\begin{proposition}\label{prop:double-duality}
The dg functors $\Hom_{\E}(P_{\bte},-)$ and $\Hom_{\E}(-,I_{\bte})$ induce quasi-isomorphisms of dg algebras
\begin{align}
\beta:R\lto\End_{\dgmod(A)}(M)^{\op} \quad\text{and}\quad \alpha:A\lto\End_{\dgmod(R)}(M)^{\op},\notag
\end{align}
respectively.
\end{proposition}
\begin{proof}
The dg functor $\Hom_{\E}(P_{\bte},-)$ induces a morphism of dg algebras
\begin{align}
\End_{\E}(I_{\bte})\lto\End_{\dgmod(A)}(M)\quad(b\longmapsto b\circ -).\notag
\end{align}
By Lemma \ref{lem:Yoneda}(3), this is a quasi-isomorphism.
Taking dg opposites, we obtain the quasi-isomorphism $\beta$, which is given by $\beta(r)=r\circ-$ for $r\in R$.
Similarly, the contravariant dg functor $\Hom_{\E}(-,I_{\bte})$ induces the morphism of dg algebras $\alpha$.
Explicitly, for homogeneous $a\in A$ and $m\in M$, the map $\alpha$ is given by $\alpha(a)(m)=(-1)^{|a||m|}m\circ a$.
The Koszul sign convention ensures that $\alpha$ is a morphism of dg algebras.
By Lemma \ref{lem:Yoneda_op}(3), $\alpha$ is a quasi-isomorphism.
\end{proof}

\begin{corollary}\label{cor:double-duality}
The following statements hold.
\begin{enumerate}
\item $R$ is quasi-isomorphic to $\REnd_{A}(T)^{\op}$ as dg algebras.
\item $A$ is quasi-isomorphic to $\REnd_{R}(T')^{\op}$ as dg algebras.
\item $A$ and $R^{\op}$ are derived Morita equivalent.
\end{enumerate}
\end{corollary}
\begin{proof}
(1) As noted before Proposition \ref{prop:double-duality}, $\REnd_{A}(T)$ can be computed by $\End_{\dgmod(A)}(M)$.
Hence the assertion follows from Proposition \ref{prop:double-duality}.

(2) Similarly, $\REnd_{R}(T')$ can be computed by $\End_{\dgmod(R)}(M)$, and hence the assertion follows from Proposition \ref{prop:double-duality}.

(3) By Proposition \ref{prop:characteristic-silting}, $T$ is a silting object in $\per A$, and hence a compact generator of $\D(A)$.
Thus $A$ and $\REnd_{A}(T)$ are derived Morita equivalent by \cite{K94}.
Since $\REnd_{A}(T)$ is quasi-isomorphic to $R^{\op}$ by (1), $A$ and $R^{\op}$ are derived Morita equivalent.
\end{proof}

We establish Ringel duality in algebraic triangulated categories.

\begin{theorem}\label{thm:Ringel_duality}
The following statements hold.
\begin{enumerate}
\item There are natural isomorphisms
\begin{align}
\RHom_{A}(F(-),T)\cong G \quad\text{and}\quad \RHom_{R}(G(-),T')\cong F\notag
\end{align}
of functors on $\thick_{\T}(\bte)$, where $\RHom_{A}(X,T)$ and $\RHom_{R}(Y,T')$ are regarded as objects of $\D(R)$ and $\D(A)$ via $\beta$ and $\alpha$, respectively.
\item The functors
\begin{align}
\RHom_{A}(-,T):\per A\lto\per R \quad\text{and}\quad \RHom_{R}(-,T'):\per R\lto\per A\notag
\end{align}
are mutually quasi-inverse contravariant triangle equivalences.
\item Up to isomorphism, $\RHom_{A}(-,T)$ and $\RHom_{R}(-,T')$ send the following objects to each other:
\begin{align}
\bde(i)\longleftrightarrow\bde'(i),\quad P_{\bde}(i)\longleftrightarrow I_{\bde'}(i),\quad I_{\bde}(i)\longleftrightarrow P_{\bde'}(i).\notag
\end{align}
\end{enumerate}
In other words, the following diagram commutes up to natural isomorphism:
\[
\begin{tikzcd}
&\thick_{\T}(\bte)\ar[ld,"F"']\ar[rd,"G"]&\\
\per A \ar[rr, shift left,"{\RHom_{A}(-,T)}"]&& \per R \ar[ll,shift left, "{\RHom_{R}(-,T')}"]
\end{tikzcd}
\]
\end{theorem}
\begin{proof}
(1) Let $X\in\thick_{\T}(\bte)$.
Consider the morphism
\begin{align}
\Hom_{\E}(X,I_{\bte})\lto\Hom_{\dgmod(A)}(\Hom_{\E}(P_{\bte},X),M)\notag
\end{align}
of right dg $R$-modules induced by $\Hom_{\E}(P_{\bte},-)$, where $R$ acts on $\Hom_{\dgmod(A)}(\Hom_{\E}(P_{\bte},X),M)$ via $\beta$.
Since $\thick_{\T}(P_{\bte})=\thick_{\T}(\bte)$, this is a quasi-isomorphism by Lemma \ref{lem:Yoneda}(3).
By definition, $\Hom_{\E}(X,I_{\bte})$ represents $GX$.
Moreover, by Lemma \ref{lem:Yoneda}(1), $\Hom_{\E}(P_{\bte},X)$ is a homotopically projective dg $A$-module representing $FX$, and hence $\Hom_{\dgmod(A)}(\Hom_{\E}(P_{\bte},X),M)$ represents $\RHom_{A}(FX,T)$.
Thus the above quasi-isomorphism induces an isomorphism $GX\cong\RHom_{A}(FX,T)$ in $\D(R)$.
Since the above quasi-isomorphism is natural in $X$, we obtain the first natural isomorphism.
Similarly, we obtain the second natural isomorphism.

(2) By (1) and Propositions \ref{prop:fullyfaithful}(1) and \ref{prop:dualfullyfaithful}(1), the functors $\RHom_{A}(-,T)$ and $\RHom_{R}(-,T')$ restrict to functors $\per A\to\per R$ and $\per R\to\per A$, which are isomorphic to $G\circ F^{-1}$ and $F\circ G^{-1}$, respectively, where $F^{-1}$ and $G^{-1}$ are quasi-inverses of $F$ and $G$.
Since $G\circ F^{-1}$ and $F\circ G^{-1}$ are mutually quasi-inverse contravariant triangle equivalences, the assertion follows.

(3) By (1), for $X\in\{\bte(i),P_{\bte}(i),I_{\bte}(i)\}$, we have $\RHom_{A}(FX,T)\cong GX$ and $\RHom_{R}(GX,T')\cong FX$.
Since $F\bte(i)=\bde(i)$ and $G\bte(i)=\bde'(i)$, the assertion follows from Propositions \ref{prop:characteristic-silting} and \ref{prop:dual-silting}.
\end{proof}

In the above setting, we say that $(A,\bde)$ and $(R,\bde')$ are \emph{Ringel dual} to each other.

\section{Derived standardly stratified algebras}

Based on the Dlab--Ringel standardization method developed in Section \ref{section:DR_std_meth}, we introduce in this section the notion of derived standardly stratified algebras as a dg analogue of standardly stratified algebras.
We study their basic properties, establish Ringel duality of derived standardly stratified algebras, and characterize their smoothness.

\subsection{Definition and basic properties}

In this subsection, we introduce derived standardly stratified algebras and give their basic properties.
We realize basic proper derived standardly stratified algebras by the standardization method in Section 2.
We also characterize derived standardly stratified algebras in terms of evaluation triangles, and establish Ringel duality for basic proper derived standardly stratified algebras.

\begin{definition}\label{def:dss}
Let $(\Lambda,\leq)$ be a finite totally ordered set.
Fix a connective dg algebra $A$ and a family $\PP=(P(i))_{i\in\Lambda}$ of objects in $\D(A)$ satisfying $\add\bigoplus_{i\in\Lambda}P(i)=\add A$.
The dg algebra $A$ is called a \emph{derived standardly stratified algebra} with respect to $(\PP,\leq)$ if there exist triangles in $\D(A)$
\begin{align}
K(i)\to P(i)\to \bde(i)\to K(i)[1] \quad(i\in\Lambda)\notag
\end{align}
satisfying the following conditions:
\begin{enumerate}
\item[(DS1)] $\End_{\D(A)}(\bde(i))$ is local for all $i\in\Lambda$.
\item[(DS2)] $\Hom_{\D(A)}(\bde(i),\bde(j)[r])=0$ for all $i>j$ and $r\in\ZZ$.
\item[(DS3)] $K(i)\in\F^{\geq 0}(\bde(>i))$ for all $i\in\Lambda$, where $\bde(>i):=\{\bde(j)\mid j>i\}$.
\end{enumerate}
We call such triangles \emph{defining triangles} of $A$.
\end{definition}

When the family $\PP$ is clear from the context, we write ``with respect to $\leq$'' in place of ``with respect to $(\PP,\leq)$'', and omit this phrase when the order is also clear.

Every standardly stratified algebra, regarded as a dg algebra concentrated in degree zero, is a derived standardly stratified algebra with respect to the same order in a natural way: the canonical short exact sequences $0\to K(i)\to P(i)\to\Delta(i)\to 0$, where $\Delta(i)$ are the standard modules, give defining triangles.
On the other hand, even if $A$ is derived standardly stratified, $H^{0}(A)$ need not be standardly stratified, as the following example illustrates.

\begin{example}\label{ex:two_cycle}
Let $Q$ be the quiver
\[
\begin{tikzcd}
1 \arrow[r, shift left, "a"] & 2 \arrow[l, shift left, "b"]
\end{tikzcd}
\]
and let $A=\kk Q/\langle ab,ba\rangle$, regarded as a dg algebra concentrated in degree zero.
Consider the order $1<2$.
Let $P(i):=e_{i}A$ be the indecomposable projective module associated with the vertex $i$.
Then we can easily check that the triangles
\begin{align}
0\to P(2)\to \bde(2)\to 0[1] \quad\text{and}\quad P(2)\xrightarrow{a\cdot -}P(1)\to\bde(1)\to P(2)[1]\notag
\end{align}
are defining triangles.
Thus $A$ is derived standardly stratified.
On the other hand, $H^{0}(A)=A$ is not standardly stratified with respect to $1<2$, since the trace of $P(2)$ in $P(1)$ is isomorphic to the simple module $S(2)$, which does not belong to $\F(P(2))$.
By symmetry, the same holds for the order $2<1$.
\end{example}

Let $A$ be a derived standardly stratified algebra with respect to $(\PP,\leq)$, with defining triangles as in Definition \ref{def:dss}, and put $\bde:=\{\bde(i)\mid i\in\Lambda\}$.
We collect basic properties of derived standardly stratified algebras.

\begin{lemma}\label{lem:thick=per}
We have $\F^{\geq 0}(\PP)=\F^{\geq 0}(\bde)$ and $\per A=\thick_{\D(A)}(\bde)$.
In particular, $\bde(i)\in\per A$ for all $i\in\Lambda$.
\end{lemma}
\begin{proof}
By (DS3) and the defining triangles, we have $P(i)\in\F^{\geq 0}(\bde)$ for all $i\in\Lambda$.
Since $\F^{\geq 0}(\bde)$ is closed under positive shifts, we obtain $\F^{\geq 0}(\PP)\subseteq\F^{\geq 0}(\bde)$.
Conversely, we show that $\bde(i)\in\F^{\geq 0}(\PP)$ by descending induction on $i$.
If $i$ is the maximal element of $\Lambda$, then $K(i)=0$, and hence $\bde(i)\cong P(i)\in\F^{\geq 0}(\PP)$.
Assume that $\bde(j)\in\F^{\geq 0}(\PP)$ for all $j>i$.
Then (DS3) gives $K(i)\in\F^{\geq 0}(\PP)$.
Since $\F^{\geq 0}(\PP)$ is closed under extensions and positive shifts, the defining triangle implies $\bde(i)\in\F^{\geq 0}(\PP)$.
Thus $\F^{\geq 0}(\bde)=\F^{\geq 0}(\PP)$.
Taking thick closures and using $\add\bigoplus_{i\in\Lambda}P(i)=\add A$, we obtain
\begin{align}
\thick_{\D(A)}(\bde) =\thick_{\D(A)}(\PP) =\thick_{\D(A)}(A) =\per A.\notag
\end{align}
This completes the proof.
\end{proof}

As shown below, the set $\bde$ is a higher standardizable set.

\begin{lemma}\label{lem:semiorth}
The following statements hold.
\begin{enumerate}
\item $\Hom_{\D(A)}(P(i),Y[r])=0$ for all $i\in\Lambda$, $Y\in\F^{\geq 0}(\bde)$ and $r>0$. In particular, the defining triangles are as in Proposition \ref{prop:const_psilt}.
\item $\Hom_{\D(A)}(\bde(i),\bde(i)[r])=0$ for all $i\in\Lambda$ and $r>0$. In particular, $\bde$ is a full higher standardizable set in $\per A$ with respect to $\leq$.
\end{enumerate}
\end{lemma}
\begin{proof}
(1) Let $Y\in\F^{\geq 0}(\bde)$.
By Lemma \ref{lem:thick=per}, we have $Y\in\F^{\geq 0}(\PP)$.
Since $A$ is connective and $P(i)\in\add A$, we obtain $\Hom_{\D(A)}(P(i),Y[r])=0$ for all $r>0$.
The last assertion follows from this and (DS3).

(2) By (1), we have $\Hom_{\D(A)}(P(i),\bde(i)[r])=0$ for all $r>0$.
On the other hand, (DS2) and (DS3) imply $\Hom_{\D(A)}(K(i),\bde(i)[r])=0$ for all $r\in\ZZ$.
Applying $\Hom_{\D(A)}(-,\bde(i))$ to the defining triangle, we obtain
\begin{align}
\Hom_{\D(A)}(\bde(i),\bde(i)[r]) \cong\Hom_{\D(A)}(P(i),\bde(i)[r])=0\notag
\end{align}
for all $r>0$.
Hence, by (DS1), (DS2) and Lemma \ref{lem:thick=per}, $\bde$ is a full higher standardizable set in $\per A$.
\end{proof}

We discuss the right minimal versions of the defining triangles and the indecomposability of the objects $P(i)$.

\begin{remark}\label{rem:dss_indec}
Assume that $A$ is proper.
We may assume that each $P(i)$ is indecomposable and that $P:=\bigoplus_{i\in\Lambda}P(i)$ is basic.
Indeed, by Remark \ref{rem:DR-minimal} and Lemma \ref{lem:semiorth}, we can replace the maps $P(i)\to\bde(i)$ by their right minimal versions, keeping the same notation.
Here, $\add P=\add A$ follows from \cite[Corollary 2.28]{AI12}, since $P$ is a silting object in $\per A$ with $\add P\subseteq\add A$.
In particular, $|\Lambda|$ equals the number of isomorphism classes of indecomposable direct summands of $A$.
\end{remark}

We also note that the class of derived standardly stratified algebras is closed under quasi-isomorphisms of dg algebras.

\begin{remark}\label{rem:dss_invariance}
Let $A$ and $B$ be dg algebras, and assume that there exists a triangle equivalence $\Psi:\per A\to\per B$ with $\Psi(\add A)=\add B$.
Fix a finite totally ordered set $(\Lambda,\leq)$ and a family $\PP=(P(i))_{i\in\Lambda}$ satisfying $\add\bigoplus_{i\in\Lambda}P(i)=\add A$, and put $\Psi\PP:=(\Psi P(i))_{i\in\Lambda}$.
Then $A$ is derived standardly stratified with respect to $(\PP,\leq)$ if and only if $B$ is derived standardly stratified with respect to $(\Psi\PP,\leq)$, since connectivity and the defining conditions are preserved under $\Psi$; note that the defining triangles lie in $\per A$ by Lemma \ref{lem:thick=per}.
In particular, this property is invariant under quasi-isomorphisms of dg algebras.
In contrast, it is not invariant under derived Morita equivalences (see Example \ref{ex:rickard}).
\end{remark}

We explain the relation with the construction in Section \ref{section:DR_std_meth}.

\begin{theorem}\label{thm:standardization}
Let $\E$ be a $\kk$-linear pretriangulated dg category, and let $\bte$ be a higher standardizable set in $\T=\H^{0}(\E)$ with respect to $\leq$ satisfying (HF).
Then the dg algebra $A=\End_{\E}(P_{\bte})$ constructed in Subsection \ref{subsection:Ringel_alg} is a proper derived standardly stratified algebra with respect to the family $(FP_{\bte}(i))_{i\in\Lambda}$ and the order $\leq$, and the triangles obtained by applying $F$ to \eqref{eq:P_tri} are defining triangles.
Conversely, every basic proper derived standardly stratified algebra with respect to $(\PP,\leq)$ is quasi-isomorphic to a dg algebra obtained in this way.
\end{theorem}
\begin{proof}
As noted in Subsection \ref{subsection:Ringel_alg}, the dg algebra $A$ is connective and proper.
We have $\bigoplus_{i\in\Lambda}FP_{\bte}(i)=FP_{\bte}\cong A$.
Applying the triangle functor $F$ to \eqref{eq:P_tri}, we obtain triangles
\begin{align}
FK_{\bte}(i)\to FP_{\bte}(i)\to\bde(i)\to FK_{\bte}(i)[1]\quad(i\in\Lambda)\notag
\end{align}
with $FK_{\bte}(i)\in\F^{\geq 0}(\bde(>i))$. Thus we have (DS3).
Conditions (DS1) and (DS2) follow from Proposition \ref{prop:fullyfaithful}(3).

Conversely, let $A$ be a basic proper derived standardly stratified algebra with respect to $(\PP,\leq)$ with defining triangles $K(i)\to P(i)\to\bde(i)\to K(i)[1]$ $(i\in\Lambda)$, and let $\E_{A}$ be the dg category of homotopically projective dg $A$-modules belonging to $\per A$, so that $\H^{0}(\E_{A})$ is triangle equivalent to $\per A$ by \cite[4.1]{K94}.
By Lemma \ref{lem:semiorth}(2), $\bde$ is a full higher standardizable set in $\per A$, and it satisfies (HF) since $A$ is proper.
By Remark \ref{rem:dss_indec}, we may assume that each $P(i)$ is indecomposable, and then $\bigoplus_{i\in\Lambda}P(i)\cong A$ since $A$ is basic.
By Lemma \ref{lem:semiorth}(1), the defining triangles can be taken as \eqref{eq:P_tri} in the construction of Subsection \ref{subsection:Ringel_alg} applied to $\E_{A}$ and $\bde$ in place of $\E$ and $\bte$.
Since $P_{\bde}=\bigoplus_{i\in\Lambda}P(i)\cong A$ in $\H^{0}(\E_{A})$, the resulting dg algebra $\End_{\E_{A}}(P_{\bde})$ is quasi-isomorphic to $\End_{\E_{A}}(A)=A$.
By Remark \ref{rem:DR-minimal}, the object $P_{\bde}$ does not depend on the choice of \eqref{eq:P_tri} up to isomorphism.
\end{proof}

\begin{corollary}\label{cor:tri_char}
Let $\T$ be an algebraic triangulated category.
Then the following statements are equivalent.
\begin{enumerate}
\item $\T$ is triangle equivalent to $\per A$ for some proper derived standardly stratified algebra $A$.
\item $\T$ admits a full higher standardizable set satisfying (HF).
\end{enumerate}
\end{corollary}
\begin{proof}
(1)$\Rightarrow$(2): Choose defining triangles of $A$ and let $\bde$ be the corresponding set. 
By Lemma \ref{lem:semiorth}(2), $\bde$ is a full higher standardizable set in $\per A$, and it satisfies (HF) since $A$ is proper.
Hence its image under a triangle equivalence $\per A\to\T$ is a full higher standardizable set in $\T$ satisfying (HF).

(2)$\Rightarrow$(1): Since $\T$ is algebraic, we may assume that $\T=\H^{0}(\E)$ for a $\kk$-linear pretriangulated dg category $\E$.
Let $\bte$ be a full higher standardizable set in $\T$ satisfying (HF).
By Theorem \ref{thm:standardization}, the dg algebra $A=\End_{\E}(P_{\bte})$ is a proper derived standardly stratified algebra.
Since $\bte$ is full, Proposition \ref{prop:fullyfaithful}(1) gives a triangle equivalence $\T=\thick_{\T}(\bte)\to\per A$.
\end{proof}

We give an explicit description of the objects $\bde(i)$ appearing in Definition \ref{def:dss}.
Let $A$ be an arbitrary dg algebra.
Let $(\Lambda,\leq)$ be a finite totally ordered set and let $\PP=(P(i))_{i\in\Lambda}$ be a family of objects in $\D(A)$ satisfying
\begin{align}
\add\bigoplus_{i\in\Lambda}P(i)=\add A.\notag
\end{align}
For each $i\in\Lambda$, put $Q_{i}:=\bigoplus_{j>i}P(j)$.
Note that $Q_{i}=0$ if $i$ is the maximal element of $\Lambda$.
Fix a homotopically projective representative of $Q_{i}$, which we denote by the same symbol, and let $E_{i}:=\End_{\dgmod(A)}(Q_{i})$.
Since $Q_{i}$ is naturally a dg $E_{i}$-$A$-bimodule, we have the derived tensor--Hom adjunction
\[
\begin{tikzcd}
\D(E_{i})\ar[rr, shift left, "{-\otimes_{E_{i}}^{\LL}Q_{i}}"]&&\D(A)\ar[ll, shift left, "{\RHom_{A}(Q_{i},-)}"]
\end{tikzcd}.
\]

Let $\S_{i}:=\Loc_{\D(A)}(Q_{i})$ be the smallest localizing subcategory of $\D(A)$ containing $Q_{i}$, and define its right perpendicular subcategory by
\begin{align}
\S_{i}^{\perp}:=\{X\in\D(A)\mid \Hom_{\D(A)}(S,X)=0
\text{ for all }S\in\S_{i}\}.\notag
\end{align}
Since $Q_{i}$ is compact, $(\S_{i},\S_{i}^{\perp})$ is a (stable) $t$-structure on $\D(A)$ (see \cite[Theorem 5.6.1]{K10} and \cite[Proposition 2.4]{M03}).

We recall the following well-known facts.

\begin{lemma}[{\cite[4.2 and 5.3]{K94}}]\label{lem:tensor_equiv}
The following statements hold.
\begin{enumerate}
\item The functor $-\otimes_{E_{i}}^{\LL}Q_{i}$ induces a triangle equivalence $\D(E_{i})\to\S_{i}$.
\item The equivalence in (1) restricts to a triangle equivalence $\per E_{i}\to\thick_{\D(A)}(Q_{i})$.
\item An object of $\S_{i}$ is compact in $\D(A)$ if and only if it belongs to $\thick_{\D(A)}(Q_{i})$.
\end{enumerate}
\end{lemma}

Let $L_{i}:=\RHom_{A}(Q_{i},P(i))\otimes_{E_{i}}^{\LL}Q_{i}$ and denote by $\ev_{i}:L_{i}\to P(i)$ the counit of the above adjunction at $P(i)$.
We call the triangle
\begin{align}
L_{i}\xrightarrow{\ev_{i}}P(i) \to\Co(\ev_{i})\to L_{i}[1]\notag
\end{align}
the \emph{evaluation triangle} of $P(i)$, and $\Co(\ev_{i})$ the \emph{derived standard module} at $i$ with respect to $(\PP,\leq)$.
If $i$ is the maximal element of $\Lambda$, then $L_{i}=0$, and hence $\Co(\ev_{i})\cong P(i)$.

We describe evaluation triangles in terms of the $t$-structure $(\S_{i},\S_{i}^{\perp})$.

\begin{lemma}\label{lem:ev_canonical}
The morphism $\RHom_{A}(Q_{i},\ev_{i})$ is an isomorphism, and we have $L_{i}\in\S_{i}$ and $\Co(\ev_{i})\in\S_{i}^{\perp}$.
In particular, the evaluation triangle of $P(i)$ is the canonical triangle of $P(i)$ with respect to $(\S_{i},\S_{i}^{\perp})$.
\end{lemma}
\begin{proof}
Since the functor $-\otimes_{E_{i}}^{\LL}Q_{i}$ is exact, preserves coproducts, and sends $E_{i}$ to $Q_{i}$, we have $L_{i}\in\S_{i}$.
By Lemma \ref{lem:tensor_equiv}(1), the functor $-\otimes_{E_{i}}^{\LL}Q_{i}$ is fully faithful.
Since $\ev_{i}$ is the counit of the adjunction, $\RHom_{A}(Q_{i},\ev_{i})$ is an isomorphism.
Consequently, we have $\RHom_{A}(Q_{i},\Co(\ev_{i}))=0$, and hence $\Co(\ev_{i})\in\S_{i}^{\perp}$.
Thus the evaluation triangle is the canonical triangle of $P(i)$ with respect to $(\S_{i},\S_{i}^{\perp})$.
\end{proof}

We give basic properties of evaluation triangles for connective dg algebras.
Let $\D^{\leq 0}(A):=\{X\in\D(A)\mid H^{j}(X)=0\text{ for all }j>0\}$.

\begin{lemma}\label{lem:ev_basic}
Assume that $A$ is connective.
Then the following statements hold.
\begin{enumerate}
\item $L_{i},\Co(\ev_{i})\in\D^{\leq 0}(A)$.
\item $\Co(\ev_{i})=0$ if and only if $P(i)\in\add Q_{i}$.
\end{enumerate}
\end{lemma}
\begin{proof}
(1) Since $P(i),Q_{i}\in\add A$, the dg algebra $E_{i}$ and the dg modules $\RHom_{A}(Q_{i},P(i))$ and $Q_{i}$ are connective.
Applying the derived K\"unneth formula \cite[Theorem 7]{Y23}, after taking a nonpositive truncation of $E_{i}$, we obtain
\begin{align}
H^{p+q}(L_{i}) \cong H^{p}(\RHom_{A}(Q_{i},P(i))) \otimes_{H^{0}(E_{i})}H^{q}(Q_{i}) \quad(p,q\geq 0).
\label{eq:ev_kunneth}
\end{align}
Since $H^{p}(\RHom_{A}(Q_{i},P(i)))=0$ and $H^{q}(Q_{i})=0$ for all $p,q>0$, we have $L_{i}\in\D^{\leq 0}(A)$.
By $P(i)\in\D^{\leq 0}(A)$, the evaluation triangle implies $\Co(\ev_{i})\in\D^{\leq 0}(A)$.

(2) Assume that $P(i)\in\add Q_{i}$.
Since the counit is an isomorphism on $\add Q_{i}$ by Lemma \ref{lem:tensor_equiv}(1), the morphism $\ev_{i}$ is an isomorphism.
Hence $\Co(\ev_{i})=0$.
Conversely, assume that $\Co(\ev_{i})=0$.
Then $\ev_{i}$ is an isomorphism.
By \eqref{eq:ev_kunneth} with $p=q=0$, we have
\begin{align}
H^{0}(L_{i})\cong \Hom_{\D(A)}(Q_{i},P(i)) \otimes_{H^{0}(E_{i})}H^{0}(Q_{i}).\notag
\end{align}
By the functoriality in \cite[Theorem 7]{Y23}, this is an isomorphism of right $H^{0}(A)$-modules.
Since $H^{0}(\ev_{i})$ is an isomorphism and $H^{0}(P(i))$ is finitely generated, there exists a surjective homomorphism
\begin{align}
H^{0}(Q_{i})^{\oplus m}\to H^{0}(P(i)) \notag
\end{align}
for some $m\geq 0$.
Since $H^{0}(P(i))$ is projective, this surjection splits.
Thus $H^{0}(P(i))\in\add H^{0}(Q_{i})$.
The equivalence $\add A\simeq\proj H^{0}(A)$ induced by $H^{0}$ gives $P(i)\in\add Q_{i}$.
\end{proof}

We show that the defining triangles of derived standardly stratified algebras are isomorphic to the evaluation triangles.

\begin{proposition}\label{prop:dss_ev}
Let $A$ be a derived standardly stratified algebra with respect to $(\PP,\leq)$, with defining triangles $K(i)\to P(i)\to\bde(i)\to K(i)[1]$.
Then, for each $i\in\Lambda$, this triangle is isomorphic to the evaluation triangle of $P(i)$.
In particular, $K(i)\cong L_{i}$ and $\bde(i)\cong\Co(\ev_{i})$.
\end{proposition}
\begin{proof}
Fix $i\in\Lambda$.
By descending induction using (DS3) of Definition \ref{def:dss}, we have
\begin{align}
\thick_{\D(A)}(Q_{i})=\thick_{\D(A)}(\bde(>i)).\notag
\end{align}
Thus $K(i)\in\S_{i}$.
Moreover, (DS2) implies $\bde(i)\in\S_{i}^{\perp}$.
On the other hand, Lemma \ref{lem:ev_canonical} gives $L_{i}\in\S_{i}$ and $\Co(\ev_{i})\in\S_{i}^{\perp}$.
Hence both the defining triangle and the evaluation triangle of $P(i)$ are canonical triangles with respect to the $t$-structure $(\S_{i},\S_{i}^{\perp})$.
By the uniqueness of such triangles, they are isomorphic.
In particular, $K(i)\cong L_{i}$ and $\bde(i)\cong\Co(\ev_{i})$.
\end{proof}

By Proposition \ref{prop:dss_ev}, the objects $\bde(i)$ are determined up to isomorphism by $A$ and $(\PP,\leq)$.
We also call $\bde(i)$ $(i\in\Lambda)$ the \emph{derived standard modules} of $A$ with respect to $(\PP,\leq)$.

We investigate when the derived standard modules associated with a connective dg algebra form a higher standardizable set.

\begin{proposition}\label{prop:ev_hss}
Assume that $A$ is connective and that the objects $P(i)$ $(i\in\Lambda)$ are pairwise non-isomorphic and have local endomorphism algebras.
Then $\{\Co(\ev_{i})\mid i\in\Lambda\}$ is a higher standardizable set in $\D(A)$ with respect to $\leq$.
\end{proposition}
\begin{proof}
Let $C_{i}:=\Co(\ev_{i})$ for each $i\in\Lambda$.
By Lemmas \ref{lem:ev_canonical} and \ref{lem:ev_basic}(1), we have $L_{i}\in\S_{i}$, $C_{i}\in\S_{i}^{\perp}$, and $L_{i},C_{i}\in\D^{\leq 0}(A)$.

We first show (S1).
Since $\Hom_{\D(A)}(P(i),L_{i}[1])=0$ and $\Hom_{\D(A)}(L_{i},C_{i}[r])=0$ for $r\in \{0,-1\}$, the evaluation triangle induces a surjective ring homomorphism $\End_{\D(A)}(P(i))\to\End_{\D(A)}(C_{i})$.
Moreover, by Lemma \ref{lem:ev_basic}(2) and the assumptions on $P(i)$, we have $C_{i}\neq 0$.
Thus $\End_{\D(A)}(C_{i})$ is local.

Next, we show (S2).
Since $P(i)\in\add A$, we have $\Hom_{\D(A)}(P(i),C_{i}[r])=0$ for all $r>0$.
Applying $\Hom_{\D(A)}(-,C_{i})$ to the evaluation triangle and using $C_{i}\in\S_{i}^{\perp}$, we obtain
\begin{align}
\Hom_{\D(A)}(C_{i},C_{i}[r]) \cong\Hom_{\D(A)}(P(i),C_{i}[r])=0\notag
\end{align}
for all $r>0$.

Finally, we show (S3).
If $i>j$, then $P(i)\in\add Q_{j}$ and $L_{i}\in \S_{i}\subseteq\S_{j}$.
Since $\S_{j}$ is closed under cones, we have $C_{i}\in \S_{j}$.
Hence $\Hom_{\D(A)}(C_{i},C_{j}[r])=0$ for all $r\in \ZZ$.
This completes the proof.
\end{proof}

Combining these results, we obtain the following characterization of derived standardly stratified algebras in terms of evaluation triangles.

\begin{theorem}\label{thm:dss_characterization}
Let $A$ be a connective dg algebra, and let $\PP=(P(i))_{i\in\Lambda}$ be a family of pairwise non-isomorphic objects with local endomorphism algebras satisfying $\add\bigoplus_{i\in\Lambda}P(i)=\add A$.
Then the following statements are equivalent.
\begin{enumerate}
\item $A$ is derived standardly stratified with respect to $(\PP,\leq)$.
\item $L_{i}\in\F^{\geq 0}(\{\Co(\ev_{j})\mid j>i\})$ for all $i\in\Lambda$.
\item $\Co(\ev_{i})\in\per A$ for all $i\in\Lambda$.
\item $\RHom_{A}(Q_{i},P(i))\in\per E_{i}$ for all $i\in\Lambda$.
\end{enumerate}
In this case, the defining triangles can be chosen to be the evaluation triangles of $P(i)$.
\end{theorem}
\begin{proof}
Let $C_{i}:=\Co(\ev_{i})$ for each $i\in\Lambda$.

(1)$\Rightarrow$(2): By Proposition \ref{prop:dss_ev}, the defining triangles are isomorphic to the evaluation triangles.
In particular, $K(i)\cong L_{i}$ and $\bde(i)\cong C_{i}$ for all $i\in\Lambda$.
Thus (2) follows from Definition \ref{def:dss}.

(2)$\Rightarrow$(1): By Proposition \ref{prop:ev_hss}, the set $\{C_{i}\mid i\in\Lambda\}$ is higher standardizable.
Combining this with (2) and $\add\bigoplus_{i\in\Lambda}P(i)=\add A$, we see that the evaluation triangles satisfy all the conditions in Definition \ref{def:dss}.
Hence (1) holds, and the defining triangles can be chosen to be the evaluation triangles.

(2)$\Rightarrow$(3): We show that $C_{i}\in\F^{\geq 0}(\{P(j)\mid j\geq i\})$ for all $i\in\Lambda$ by descending induction on $i$.
If $i$ is the maximal element of $\Lambda$, then $L_{i}=0$, and hence $C_{i}\cong P(i)$.
Assume that $C_{j}\in\F^{\geq 0}(\{P(l)\mid l\geq j\})$ for all $j>i$.
Then $L_{i}\in\F^{\geq 0}(\{P(j)\mid j>i\})$ by (2).
Since $\F^{\geq 0}(\{P(j)\mid j\geq i\})$ is closed under extensions and positive shifts, the evaluation triangle of $P(i)$ implies $C_{i}\in\F^{\geq 0}(\{P(j)\mid j\geq i\})$.
In particular, $C_{i}\in\per A$.

(3)$\Leftrightarrow$(4): The assertion follows from
\begin{align}
C_{i}\in\per A
&\Leftrightarrow L_{i}\in\per A&&(\text{by }P(i)\in\per A)\notag\\
&\Leftrightarrow L_{i}\in\thick_{\D(A)}(Q_{i})&&(\text{by }L_{i}\in\S_{i}\text{ and Lemma \ref{lem:tensor_equiv}(3)})\notag\\
&\Leftrightarrow \RHom_{A}(Q_{i},P(i))\in\per E_{i}&&(\text{by Lemma \ref{lem:tensor_equiv}(1) and (2)}).\notag
\end{align}

(3)$\Rightarrow$(2): We use descending induction on $i$.
If $i$ is the maximal element of $\Lambda$, then $L_{i}=0$.
Assume that $L_{j}\in\F^{\geq 0}(\{C_{l}\mid l>j\})$ for all $j>i$.
By Proposition \ref{prop:ev_hss}, $\{C_{j}\mid j>i\}$ is a higher standardizable set, and hence $\F^{\geq 0}(\{C_{j}\mid j>i\})$ is closed under direct summands.
For $j>i$, the evaluation triangle of $P(j)$ gives $P(j)\in\F^{\geq 0}(\{C_{l}\mid l\geq j\})$.
Hence $\F^{\geq 0}(\{P(j)\mid j>i\})\subseteq\F^{\geq 0}(\{C_{j}\mid j>i\})$.
On the other hand, (3) gives $L_{i}\in\thick_{\D(A)}(Q_{i})$, as in the proof of (3)$\Leftrightarrow$(4).
Since $Q_{i}=\bigoplus_{j>i}P(j)$ is presilting, it follows from \cite[Proposition 2.8(a)]{IY18} that there exists a triangle
\begin{align}
L'_{i}\to L_{i}\to L''_{i}\to L'_{i}[1]\notag
\end{align}
such that $L'_{i}\in \F^{\leq 0}(\{P(j)\mid j>i\})[-1]$ and $L''_{i}\in \F^{\geq 0}(\{P(j)\mid j>i\})$.
Since $\RHom_{A}(Q_{i},\ev_{i})$ is an isomorphism by Lemma \ref{lem:ev_canonical} and $A$ is connective, we have
\begin{align}
\Hom_{\D(A)}(Q_{i},L_{i}[r])\cong\Hom_{\D(A)}(Q_{i},P(i)[r])=0\notag
\end{align}
for all $r>0$.
Hence $L_{i}$ is a direct summand of $L''_{i}$.
Thus we obtain $L_{i}\in \F^{\geq 0}(\{C_{j}\mid j>i\})$.
\end{proof}

We apply Theorem \ref{thm:dss_characterization} to finite-dimensional algebras.
Let $A$ be a basic finite-dimensional algebra, regarded as a dg algebra concentrated in degree zero, let $\{e_{i}\mid i\in\Lambda\}$ be a complete set of primitive orthogonal idempotents of $A$, and put $P(i):=e_{i}A$ and $\e_{i}:=\sum_{j>i}e_{j}$. Note that the objects $P(i)$ are pairwise non-isomorphic and have local endomorphism algebras, and $\add\bigoplus_{i\in\Lambda}P(i)=\add A$.

\begin{corollary}\label{cor:fd_criterion}
The algebra $A$ is derived standardly stratified with respect to the family $(e_{i}A)_{i\in\Lambda}$ and the order $\leq$ if and only if the right $\e_{i}A\e_{i}$-module $e_{i}A\e_{i}$ has finite projective dimension for every $i\in\Lambda$.
\end{corollary}
\begin{proof}
Since $Q_{i}=\e_{i}A$ is a projective $A$-module, it is homotopically projective. 
Thus $E_{i}=\End_{A}(\e_{i}A)=\e_{i}A\e_{i}$ and $\RHom_{A}(Q_{i},P(i))=\Hom_{A}(\e_{i}A,e_{i}A)=e_{i}A\e_{i}$, concentrated in degree zero.
It is known that a finitely generated module over a finite-dimensional algebra is perfect if and only if it has finite projective dimension.
Thus $\RHom_{A}(Q_{i},P(i))\in \per E_{i}$ if and only if $e_{i}A\e_{i}$ has finite projective dimension.
Hence the assertion follows from the equivalence (1)$\Leftrightarrow$(4) in Theorem \ref{thm:dss_characterization}.
\end{proof}

\begin{example}
Let $n\geq 2$, and let $Q$ be the cyclic quiver
\[
\begin{tikzcd}
1 \arrow[r, "a_{1}"] & 2 \arrow[r, "a_{2}"] & \cdots \arrow[r, "a_{n-1}"] & n \arrow[lll, bend left=20, "a_{n}"]
\end{tikzcd}.
\]
Let $A=\kk Q/J^{2}$, where $J$ is the arrow ideal of $\kk Q$, and let $e_{i}$ be the idempotent corresponding to the vertex $i$.
Then $A$ is derived standardly stratified with respect to the family $(e_{i}A)_{i=1}^{n}$ and every total order $\leq$ on $\{1,\ldots,n\}$.
Indeed, for each vertex $i$, the algebra $\e_{i}A\e_{i}$ is the radical square zero algebra of the full subquiver of $Q$ on the vertices $j>i$.
This subquiver has no oriented cycles since it does not contain $i$.
Hence $\e_{i}A\e_{i}$ has finite global dimension, and the assertion follows from Corollary \ref{cor:fd_criterion}.
On the other hand, $\gldim A=\infty$.
\end{example}

We establish Ringel duality for basic proper derived standardly stratified algebras.
Let $A$ be a basic proper derived standardly stratified algebra with respect to $(\PP,\leq)$, with defining triangles $K(i)\to P(i)\to\bde(i)\to K(i)[1]$ $(i\in\Lambda)$ such that each $P(i)$ is indecomposable (see Remark \ref{rem:dss_indec}).
Let $\E_{A}$ be the dg category of homotopically projective dg $A$-modules belonging to $\per A$, so that $\H^{0}(\E_{A})$ is triangle equivalent to $\per A$ by \cite[4.1]{K94}.
By Lemma \ref{lem:semiorth}, we may apply Subsection \ref{subsection:Ringel_alg} to $\E_{A}$ and $\bde$ in place of $\E$ and $\bte$, taking the defining triangles as \eqref{eq:P_tri} and $P_{\bde}=A$, so that $F\cong\id$.
Let $R$, $G$ and $\bde'$ be as in Subsection \ref{subsection:Ringel_alg}.
Put $\PP':=(GI_{\bde}(i))_{i\in\Lambda}$.

\begin{theorem}\label{thm:Ringel_dss}
The following statements hold.
\begin{enumerate}
\item $R$ is a basic proper derived standardly stratified algebra with respect to $(\PP',\leq^{\op})$.
The triangles obtained by applying $G$ to \eqref{eq:I_tri} are defining triangles of $R$.
\item $(A,\bde)$ and $(R,\bde')$ are Ringel dual to each other.
\end{enumerate}
\end{theorem}
\begin{proof}
(1) As noted in Subsection \ref{subsection:Ringel_alg}, the dg algebra $R$ is connective and proper.
We obtain $\bigoplus_{i\in\Lambda}GI_{\bde}(i)=GI_{\bde}\cong R$.
Since $\bde$ is a higher standardizable set satisfying (HF), we have triangles
\begin{align}
\bde(i)\to I_{\bde}(i)\to C_{\bde}(i)\to\bde(i)[1]\quad(i\in\Lambda)\notag
\end{align}
as in \eqref{eq:I_tri}, where $C_{\bde}(i)\in\F^{\leq 0}(\bde(<i))$.
Applying $G$ to these triangles, we obtain triangles
\begin{align}
GC_{\bde}(i)\to GI_{\bde}(i)\to\bde'(i)\to GC_{\bde}(i)[1]\quad(i\in\Lambda).\notag
\end{align}
Since $G$ is a contravariant triangle functor, we have $GC_{\bde}(i)\in\F^{\geq 0}(\bde'(>^{\op}i))$, and hence (DS3) holds.
Conditions (DS1) and (DS2) with respect to $\leq^{\op}$ follow from Proposition \ref{prop:dualfullyfaithful}(3).
Hence $R$ is a proper derived standardly stratified algebra with respect to $(\PP',\leq^{\op})$, with these triangles as defining triangles.
Moreover, $R$ is basic, since $R\cong P_{\bde'}$ by Proposition \ref{prop:dual-silting}.

(2) Since $P_{\bde}=A$ and $F\cong\id$, this holds by definition.
\end{proof}

\subsection{Smoothness}

In this subsection, we characterize smoothness of proper derived standardly stratified algebras, and discuss its consequences for full exceptional sequences and derived Morita equivalences.
For the characterization, we recall the notions of projective dimension and global dimension for connective dg algebras (\cite{Y13,M21}).
Let $A$ be a connective dg algebra.
We denote by $\Mod H^{0}(A)$ the category of right $H^{0}(A)$-modules, and by $\mod H^{0}(A)$ its full subcategory consisting of finitely presented modules.
We regard $H^{0}(A)$-modules as objects of $\D(A)$.
For $M\in\mod H^{0}(A)$, we define the projective dimension of $M$ by
\begin{align}
\pd_{A}M:=\inf\{d\geq 0\mid \Hom_{\D(A)}(M,N[r])=0 \text{ for all }N\in\Mod H^{0}(A)\text{ and }r>d\}.\notag
\end{align}
Following \cite[Definition 2.5]{H26}, we define the (right) global dimension of $A$ by
\begin{align}
\gldim A:=\sup\{\pd_{A}M\mid M\in\mod H^{0}(A)\}.\notag
\end{align}

We obtain the following characterization of smoothness.

\begin{theorem}\label{thm:dgqh}
Assume that $\kk$ is algebraically closed, and let $A$ be a proper derived standardly stratified algebra with respect to $(\PP,\leq)$, with defining triangles $K(i)\to P(i)\to\bde(i)\to K(i)[1]$ $(i\in\Lambda)$.
Then the following statements are equivalent.
\begin{enumerate}
\item $\bde(i)$ is exceptional for all $i\in \Lambda$.
\item $A$ is smooth.
\item $\gldim A<\infty$.
\item $\per A$ admits a full exceptional sequence.
\end{enumerate}
\end{theorem}
\begin{proof}
Write $n=|\Lambda|$ and identify $\Lambda$ with $[1,n]$ via the order-preserving bijection.
For each $i\in[1,n]$, put $\C_{i}:=\thick_{\D(A)}(\bde(i))$ and $B_{i}:=\REnd_{A}(\bde(i))$.
Since $A$ is proper and $\bde(i)\in\per A$, each $B_{i}$ is proper.
Moreover, $B_{i}$ is connective by Lemma \ref{lem:semiorth}(2), and $H^{0}(B_{i})$ is local by (DS1).
Since the singleton $\{\bde(i)\}$ is a higher standardizable set in $\per A$ satisfying (HF), it follows from Proposition \ref{prop:fullyfaithful}(1) that there exists a triangle equivalence $\C_{i}\to\per B_{i}$ induced by a dg functor.
By (DS2) and Lemma \ref{lem:thick=per}, we have a semiorthogonal decomposition $\per A=\langle\C_{1},\ldots,\C_{n}\rangle$, that is, $\Hom_{\D(A)}(\C_{i},\C_{j})=0$ for all $i>j$, and $\per A$ is the smallest triangulated subcategory of $\D(A)$ containing $\C_{1},\ldots,\C_{n}$.
By \cite[Proposition 2.11]{E17} (see also \cite{LS14}) and the invariance of smoothness under derived Morita equivalences (see \cite[Section 4.7]{K06}), $A$ is smooth if and only if $B_{i}$ is smooth for all $i\in[1,n]$.

(1)$\Rightarrow$(2): Since $\kk$ is algebraically closed and $B_{i}$ is proper, (1) implies $H^{0}(B_{i})\cong\kk$ and $H^{r}(B_{i})=0$ for all $r\neq 0$.
Thus the canonical morphism $\kk\to B_{i}$ is a quasi-isomorphism.
Since $\kk$ is smooth, so is each $B_{i}$. Therefore $A$ is smooth.

(2)$\Rightarrow$(1): Assume that $A$ is smooth.
Then each $B_{i}$ is smooth, proper and connective with $H^{0}(B_{i})$ local.
By \cite[Corollary 7.5]{RS22}, there exists a triangle equivalence $F_{i}:\per B_{i}\to\per\kk$.
Since $\End_{\per\kk}(F_{i}(B_{i}))\cong H^{0}(B_{i})$ is local, we have $F_{i}(B_{i})\cong\kk[m_{i}]$ for some $m_{i}\in\ZZ$.
Since $B_{i}=\REnd_{A}(\bde(i))$, we have
\begin{align}
\Hom_{\D(A)}(\bde(i),\bde(i)[r])
&\cong\Hom_{\per B_{i}}(B_{i},B_{i}[r])\cong\Hom_{\per\kk}(F_{i}(B_{i}),F_{i}(B_{i})[r])\notag\\
&\cong
\begin{cases}
\kk & (r=0),\\
0 & (r\neq 0).
\end{cases}\notag
\end{align}
Thus $\bde(i)$ is exceptional.

(2)$\Leftrightarrow$(3): This follows from \cite[Propositions 2.9 and 2.10]{H26}, since $A$ is proper and $\kk$ is algebraically closed.

(1)$\Rightarrow$(4): By (1), (DS2) and Lemma \ref{lem:thick=per}, $(\bde(1),\ldots,\bde(n))$ is a full exceptional sequence in $\per A$.

(4)$\Rightarrow$(2): Let $(E_{1},\ldots,E_{m})$ be a full exceptional sequence in $\per A$.
Then the arguments in the first paragraph and in (1)$\Rightarrow$(2) work with $\bde(i)$ replaced by $E_{j}$.
Hence $A$ is smooth.
\end{proof}

In view of Theorem \ref{thm:dgqh}, the four statements suggest four possible definitions of a dg analogue of quasi-hereditary algebras.
We adopt (1) in Theorem \ref{thm:dgqh}, which is formulated in terms of the derived standard modules.
Thus, a derived standardly stratified algebra with respect to $(\PP,\leq)$ is called a \emph{derived quasi-hereditary algebra} if all the corresponding derived standard modules are exceptional.

Several results on derived standardly stratified algebras in the previous subsection have analogues for derived quasi-hereditary algebras.

\begin{remark}\label{rem:dqh}
The following statements hold.
\begin{enumerate}
\item Every quasi-hereditary algebra is derived quasi-hereditary with respect to the same order.
\item In the situation of Remark \ref{rem:dss_invariance}, $A$ is derived quasi-hereditary with respect to $(\PP,\leq)$ if and only if $B$ is derived quasi-hereditary with respect to $(\Psi\PP,\leq)$.
\item In the situation of Theorem \ref{thm:Ringel_dss}, $A$ is derived quasi-hereditary with respect to $(\PP,\leq)$ if and only if $R$ is derived quasi-hereditary with respect to $(\PP',\leq^{\op})$.
\end{enumerate}
\end{remark}

In contrast to Remarks \ref{rem:dss_invariance} and \ref{rem:dqh}(2), neither the class of derived standardly stratified algebras nor that of derived quasi-hereditary algebras is closed under derived Morita equivalences, even among finite-dimensional algebras.
The following example is based on a construction of Rickard \cite{R19}.

\begin{example}\label{ex:rickard}
Let $\kk$ be an algebraically closed field, and let $A=\kk Q/I$, where $Q$ is the quiver
\[
\begin{tikzcd}
1 \arrow[r, bend left=40, "\alpha"] & 2 \arrow[l, bend left=40, "\beta_{1}"] \arrow[l, "\beta_{2}"']
\end{tikzcd}
\]
and $I=\langle\beta_{1}\alpha,\beta_{2}\alpha\rangle$.
Put $P(i):=e_{i}A$ and $S(i):=\operatorname{top}P(i)$, where $e_{i}$ is the idempotent corresponding to the vertex $i$.
Then $A$ is quasi-hereditary with respect to the order $1<2$, with usual standard modules $\Delta(1)=S(1)$ and $\Delta(2)=P(2)$.
In particular, $A$ is derived quasi-hereditary by Remark \ref{rem:dqh}(1).
By \cite{R19}, $T=P(1)[1]\oplus X$ is a tilting object in $\per A$, where $X=(P(1)^{\oplus2}\xrightarrow{(\beta_{1},\beta_{2})}P(2))$ with $P(2)$ in degree zero.
Thus $B=\End_{\per A}(T)$ is a finite-dimensional algebra derived Morita equivalent to $A$.
Since $A$ is smooth by Theorem \ref{thm:dgqh}, so is $B$.
By a direct calculation, the endomorphism algebras of $P(1)[1]$ and $X$ are local but not division algebras.
Suppose that $B$ is derived standardly stratified with respect to some pair $(\PP_{B},\leq)$, where $\PP_{B}=(P_{B}(i))_{i\in\Lambda}$, and let $n$ be the maximal element of $\Lambda$.
Then $P_{B}(n)\cong\bde_{B}(n)$.
Since $B$ is smooth, Theorem \ref{thm:dgqh} implies that $\End_{B}(P_{B}(n))\cong\End_{\D(B)}(\bde_{B}(n))$ is a division algebra.
However, $\End_{B}(P_{B}(n))$ is isomorphic to either $\End_{\per A}(P(1)[1])$ or $\End_{\per A}(X)$, a contradiction.
Therefore $B$ is neither derived standardly stratified nor derived quasi-hereditary.
\end{example}

The following corollary is an analogue of Corollary \ref{cor:tri_char} for derived quasi-hereditary algebras.

\begin{corollary}\label{cor:exc_char}
Let $\T$ be an algebraic triangulated category.
Then the following statements are equivalent.
\begin{enumerate}
\item $\T$ is triangle equivalent to $\per A$ for some proper derived quasi-hereditary algebra $A$.
\item $\T$ admits a full exceptional sequence satisfying (HF).
\end{enumerate}
\end{corollary}
\begin{proof}
(1)$\Rightarrow$(2): Choose a pair $(\PP,\leq)$ for which $A$ is derived quasi-hereditary. By (DS2) and Lemma \ref{lem:thick=per}, the corresponding derived standard modules of $A$ form a full exceptional sequence in $\per A$.
Moreover, it satisfies (HF) since $A$ is proper.

(2)$\Rightarrow$(1): Let $\bte=(E_{1},\ldots,E_{n})$ be a full exceptional sequence in $\T$ satisfying (HF).
By Example \ref{ex:exceptional} and the proof of Corollary \ref{cor:tri_char}, $\T$ is triangle equivalent to $\per A$ for the proper derived standardly stratified algebra $A=\End_{\E}(P_{\bte})$.
By Theorem \ref{thm:standardization}, the derived standard modules of $A$ are $FE_{1},\ldots,FE_{n}$, which are exceptional.
Hence $A$ is derived quasi-hereditary.
\end{proof}

The following examples show that there exist derived quasi-hereditary algebras which are not derived Morita equivalent to any quasi-hereditary algebra.

\begin{example}\label{ex:Kalck}
Let $\kk$ be an algebraically closed field, and let $A=\kk Q/I$, where $Q$ is the quiver
\[
\begin{tikzcd}
& 2 \arrow[dl, "\alpha_{1}"'] & \\
3 \arrow[rr, shift left, "\gamma"] \arrow[rr, shift right, "\beta"'] & & 1 \arrow[ul, "\alpha_{2}"']
\end{tikzcd}
\]
and $I=\langle\beta\alpha_{2},\alpha_{1}\gamma,\alpha_{2}\alpha_{1}\rangle$.
This algebra is denoted by $A_{2}$ in \cite{Ka17}.
By \cite[Example 2.23 and Proposition 5.6]{Ka17}, $\per A$ admits a full exceptional sequence, but $A$ is not derived equivalent to any quasi-hereditary algebra.
Since $A$ is proper, this full exceptional sequence satisfies (HF).
By Corollary \ref{cor:exc_char}, $\per A$ is triangle equivalent to $\per B$ for some proper derived quasi-hereditary algebra $B$.
This triangle equivalence is induced by a dg functor (see Proposition \ref{prop:fullyfaithful}(1)), and hence $A$ and $B$ are derived Morita equivalent.
Since $A$ is not derived equivalent to any quasi-hereditary algebra, $B$ is not derived Morita equivalent to any quasi-hereditary algebra.
Moreover, this example shows that Corollary \ref{cor:exc_char} does not hold if derived quasi-hereditary algebras are replaced by quasi-hereditary algebras.
\end{example}

\begin{example}\label{ex:RS}
Let $\kk$ be an algebraically closed field, and let $A$ be the path algebra of the quiver
\[
\begin{tikzcd}
1 \arrow[r, shift left, "a"] \arrow[r, shift right, "b"'] & 2
\end{tikzcd}
\]
with $|a|=0$ and $|b|=-2$, regarded as a dg algebra with zero differential.
By \cite[Example 4.8]{RS22}, $A$ is proper, connective and smooth, but not derived Morita equivalent to any finite-dimensional algebra.
Put $P(i)=e_{i}A$, where $e_{i}$ is the idempotent corresponding to the vertex $i$, and consider the order $1<2$.
Since $2$ is maximal, we have $L_{2}=0$, and hence $\Co(\ev_{2})\cong P(2)$.
Since $\RHom_{A}(P(2),P(1))\cong\kk\oplus\kk[2]$, we have $L_{1}\cong\Co(\ev_{2})\oplus\Co(\ev_{2})[2]$.
Hence $A$ is derived standardly stratified by Theorem \ref{thm:dss_characterization}, and derived quasi-hereditary by Theorem \ref{thm:dgqh}.
In particular, $A$ is not derived Morita equivalent to any quasi-hereditary algebra.
\end{example}

The dg algebra in Example \ref{ex:RS} is the path algebra of an acyclic graded quiver.
More generally, we show that dg algebras arising from finite acyclic graded quivers with arrows of non-positive degrees, possibly with homogeneous relations, are derived quasi-hereditary with respect to every total order on the vertices.
For this purpose, we introduce directed decompositions.
For idempotents $e,e'\in Z^{0}(A)$ of a dg algebra $A$, we regard $eAe'$ as a complex, and note that $\Hom_{\D(A)}(e'A,eA[r])\cong H^{r}(eAe')$ for all $r\in\ZZ$.

\begin{definition}\label{def:directed}
Let $A$ be a dg algebra, and let $\{e_{i}\mid i\in\Lambda\}$ be a set of pairwise orthogonal idempotents in $Z^{0}(A)$ satisfying $\sum_{i\in\Lambda}e_{i}=1$, where $\Lambda$ is a finite set.
We call the decomposition $1=\sum_{i\in\Lambda}e_{i}$ \emph{directed} if there exists a total order $\trianglelefteq$ on $\Lambda$ satisfying the following conditions.
\begin{enumerate}
\item[(D1)] $H^{r}(e_{i}Ae_{j})=0$ for all $r\in\ZZ$ unless $i\trianglelefteq j$.
\item[(D2)] For each $i\in\Lambda$, the map $\kk\to e_{i}Ae_{i}$ ($\lambda\mapsto\lambda e_{i}$) is a quasi-isomorphism.
\end{enumerate}
\end{definition}

If $A$ admits a direct sum decomposition in $\per A$, then it does not necessarily come from idempotents in $Z^{0}(A)$.
The following remark shows that it does after replacing $A$ by a derived Morita equivalent dg algebra.

\begin{remark}\label{rem:directed}
Let $A$ be a dg algebra such that $A\cong\bigoplus_{i\in\Lambda}P(i)$ in $\per A$ for some objects $P(i)$ $(i\in\Lambda)$, where $\Lambda$ is a finite set.
Let $\til{P}(i)$ be a homotopically projective representative of $P(i)$, and put $A':=\End_{\dgmod(A)}(\bigoplus_{i\in\Lambda}\til{P}(i))$.
Then the projections $e_{i}\in Z^{0}(A')$ onto $\til{P}(i)$ are pairwise orthogonal idempotents satisfying $\sum_{i\in\Lambda}e_{i}=1$.
Since $\bigoplus_{i\in\Lambda}\til{P}(i)\cong A$ in $\D(A)$, Lemma \ref{lem:tensor_equiv}(2) applied to $\bigoplus_{i\in\Lambda}\til{P}(i)$ instead of $Q_{i}$ gives a triangle equivalence
\begin{align}
F:=-\otimes^{\LL}_{A'}\bigoplus_{i\in\Lambda}\til{P}(i)\colon\per A'\to\thick_{\D(A)}\Bigl(\bigoplus_{i\in\Lambda}\til{P}(i)\Bigr)=\per A\notag
\end{align}
with $F(e_{i}A')\cong\til{P}(i)\cong P(i)$ for all $i\in\Lambda$.
In particular, $F(\add A')=\add A$.
Put $\PP:=(P(i))_{i\in\Lambda}$ and $\PP':=(e_{i}A')_{i\in\Lambda}$.
Thus, by Remarks \ref{rem:dss_invariance} and \ref{rem:dqh}(2), for every total order $\leq$ on $\Lambda$, $A$ is derived quasi-hereditary with respect to $(\PP,\leq)$ if and only if $A'$ is derived quasi-hereditary with respect to $(\PP',\leq)$.
\end{remark}

We describe directed decompositions in terms of exceptional sequences.

\begin{lemma}\label{lem:directed}
Let $A$ be a dg algebra with a set $\{e_{i}\mid i\in\Lambda\}$ of pairwise orthogonal idempotents in $Z^{0}(A)$ satisfying $\sum_{i\in\Lambda}e_{i}=1$, and put $P(i):=e_{i}A$.
Let $\trianglelefteq$ be a total order on $\Lambda$, and write $\Lambda=\{i_{1},\ldots,i_{n}\}$ with $i_{n}\triangleleft\cdots\triangleleft i_{1}$.
Then (D1) and (D2) hold if and only if $(P(i_{1}),\ldots,P(i_{n}))$ is a full exceptional sequence in $\per A$ such that $\End_{\D(A)}(P(i))\cong\kk$ for all $i\in\Lambda$.
In this case, if $A$ is proper, then it is smooth.
\end{lemma}
\begin{proof}
By the isomorphism $\Hom_{\D(A)}(P(j),P(i)[r])\cong H^{r}(e_{i}Ae_{j})$, (D1) and (D2) hold if and only if $(P(i_{1}),\ldots,P(i_{n}))$ is an exceptional sequence in $\per A$ such that $\End_{\D(A)}(P(i))\cong\kk$ for all $i\in\Lambda$.
In this case, the exceptional sequence is full since $\bigoplus_{i\in\Lambda}P(i)=A$.

Assume that (D1) and (D2) hold and that $A$ is proper.
By (D2), $\REnd_{A}(P(i))=e_{i}Ae_{i}$ is quasi-isomorphic to $\kk$ for all $i\in\Lambda$.
Hence $A$ is smooth, since the arguments in the first paragraph and in (1)$\Rightarrow$(2) of the proof of Theorem \ref{thm:dgqh} work with $\bde(i)$ replaced by $P(i)$.
\end{proof}

We show that a proper connective dg algebra with a directed decomposition $1=\sum_{i\in\Lambda}e_{i}$ is derived quasi-hereditary with respect to the family $(e_{i}A)_{i\in\Lambda}$ and every total order on $\Lambda$.

\begin{theorem}\label{thm:directed}
Let $A$ be a proper connective dg algebra with a directed decomposition $1=\sum_{i\in\Lambda}e_{i}$, and put $P(i):=e_{i}A$ for $i\in\Lambda$ and $\PP:=(P(i))_{i\in\Lambda}$.
Then, for every total order $\leq$ on $\Lambda$, the dg algebra $A$ is derived quasi-hereditary with respect to $(\PP,\leq)$.
\end{theorem}
\begin{proof}
Fix a total order $\trianglelefteq$ on $\Lambda$ satisfying (D1) and (D2), and a total order $\leq$ on $\Lambda$.
For each $i\in\Lambda$, put $\e_{i}:=\sum_{j>i}e_{j}$.
Since $\e_{i}A$ is homotopically projective, we may take $Q_{i}=\e_{i}A$, so that $E_{i}=\End_{\dgmod(A)}(\e_{i}A)=\e_{i}A\e_{i}$ and $\RHom_{A}(Q_{i},P(i))=e_{i}A\e_{i}$.

We first show that $A$ is derived standardly stratified with respect to $(\PP,\leq)$.
By Lemma \ref{lem:directed}, the objects $P(i)$ form an exceptional sequence in a suitable order and satisfy $\End_{\D(A)}(P(i))\cong\kk$.
In particular, they are pairwise non-isomorphic and have local endomorphism algebras.
For each $i\in\Lambda$, we have $e_{j}E_{i}e_{l}=e_{j}Ae_{l}$ for all $j,l>i$, so that the decomposition $\e_{i}=\sum_{j>i}e_{j}$ of the unit of $E_{i}$ is directed.
Moreover, $E_{i}$ is proper since so is $A$.
Hence $E_{i}$ is smooth by Lemma \ref{lem:directed}.
Since $A$ is proper, the dg $E_{i}$-module $e_{i}A\e_{i}$ has finite-dimensional total cohomology.
Since $E_{i}$ is smooth, it belongs to $\per E_{i}$ by \cite[Lemma 3.8(ii)]{KS25}.
Thus the assertion follows from Theorem \ref{thm:dss_characterization}.

It remains to show that $\bde(i)$ is exceptional for each $i\in\Lambda$.
Fix $i\in\Lambda$.
Let
\begin{align}
L_{i}\xrightarrow{\ev_{i}}P(i)\to\bde(i)\to L_{i}[1]\notag
\end{align}
be the evaluation triangle of $P(i)$, which is a defining triangle by Theorem \ref{thm:dss_characterization}.
Put $Q_{i}^{+}:=\bigoplus_{j>i,\ i\triangleleft j}P(j)$ and $Q_{i}^{-}:=\bigoplus_{j>i,\ j\triangleleft i}P(j)$, so that $Q_{i}=Q_{i}^{+}\oplus Q_{i}^{-}$.
By (D1), we have
\begin{align}\label{eq:directed_vanishing}
\RHom_{A}(Q_{i}^{-},P(i))=0,\quad \RHom_{A}(Q_{i}^{-},Q_{i}^{+})=0,\quad \RHom_{A}(P(i),Q_{i}^{+})=0.
\end{align}
Put $\S_{i}:=\Loc_{\D(A)}(Q_{i})$ and $\U_{i}:=\Loc_{\D(A)}(Q_{i}^{+})$.
Recall that $(\S_{i},\S_{i}^{\perp})$ and $(\U_{i},\U_{i}^{\perp})$ are $t$-structures since $Q_{i}$ and $Q_{i}^{+}$ are compact.
Let
\begin{align}
L\to P(i)\to C\to L[1]\notag
\end{align}
be the canonical triangle of $P(i)$ with respect to $(\U_{i},\U_{i}^{\perp})$.
By the first two equalities in \eqref{eq:directed_vanishing}, we have $\RHom_{A}(Q_{i}^{-},C)=0$.
Hence $C\in\S_{i}^{\perp}$.
Since $L\in\U_{i}\subseteq\S_{i}$ and $C\in\S_{i}^{\perp}$, this triangle is the canonical triangle of $P(i)$ with respect to $(\S_{i},\S_{i}^{\perp})$.
By Lemma \ref{lem:ev_canonical} and the uniqueness of such triangles, it is isomorphic to the evaluation triangle of $P(i)$.
In particular, $L_{i}\in\U_{i}$.
By the third equality in \eqref{eq:directed_vanishing}, we have $\RHom_{A}(P(i),L_{i})=0$.
Moreover, $L_{i}\in\S_{i}$ and $\bde(i)\in\S_{i}^{\perp}$ imply $\RHom_{A}(L_{i},\bde(i))=0$.
Applying $\RHom_{A}(-,\bde(i))$ and $\RHom_{A}(P(i),-)$ to the evaluation triangle, we obtain isomorphisms in $\D(\kk)$
\begin{align}
\RHom_{A}(\bde(i),\bde(i)) \cong\RHom_{A}(P(i),\bde(i)) \cong\RHom_{A}(P(i),P(i)).\notag
\end{align}
Since $\RHom_{A}(P(i),P(i))$ is represented by $e_{i}Ae_{i}$, condition (D2) gives
\begin{align}
\Hom_{\D(A)}(\bde(i),\bde(i)[r]) \cong H^{r}(e_{i}Ae_{i}) \cong
\begin{cases}
\kk & (r=0),\\
0 & (r\neq 0)
\end{cases}\notag
\end{align}
for all $r\in\ZZ$.
Thus $\bde(i)$ is exceptional, and hence $A$ is derived quasi-hereditary with respect to $(\PP,\leq)$.
\end{proof}

We apply Theorem \ref{thm:directed} to acyclic graded quivers.
Let $Q$ be a finite acyclic quiver.
By \cite{DR89}, the path algebra $\kk Q$ is quasi-hereditary with respect to every total order on the set of vertices of $Q$, whereas $\kk Q/I$ is not for a nonzero admissible ideal $I$.
Note that $\kk Q/I$ is quasi-hereditary with respect to the total order given by any sink sequence or source sequence of the vertices of $Q$.
In contrast, $\kk Q/I$ is derived quasi-hereditary with respect to every total order on the vertices in the connective dg setting considered in the following corollary.

\begin{corollary}\label{cor:acyclic}
Let $Q$ be a finite acyclic graded quiver whose arrows have non-positive degrees, let $J$ be the arrow ideal of $\kk Q$, and let $I\subseteq J^{2}$ be a homogeneous ideal.
Let $A=\kk Q/I$, endowed with an arbitrary differential.
Then the idempotents $e_{i}$ corresponding to the vertices $i$ of $Q$ give a directed decomposition of $A$.
In particular, $A$ is derived quasi-hereditary with respect to every total order on the set of vertices of $Q$.
\end{corollary}
\begin{proof}
Since $Q$ is acyclic, there exists a total order $\trianglelefteq$ on the set of vertices of $Q$ such that $i\trianglelefteq j$ whenever there exists a path from $i$ to $j$.
Since $A$ is spanned by the images of the paths in $Q$ and all arrows have non-positive degrees, we have $A^{r}=0$ for all $r>0$.
In particular, $d(e_{i})=0$, and hence $e_{i}\in Z^{0}(A)$ and $e_{i}Ae_{j}$ is a subcomplex of $A$.
Moreover, $e_{i}Ae_{j}$ is spanned by the images of the paths from $i$ to $j$, and $e_{i}Ae_{i}=\kk e_{i}\neq 0$ since $I\subseteq J^{2}$.
Thus (D1) and (D2) hold for $\trianglelefteq$.
Since $Q$ has only finitely many paths, $A$ is finite-dimensional, and hence proper.
Moreover, $A$ is connective since $A^{r}=0$ for all $r>0$.
Hence the assertion follows from Theorem \ref{thm:directed}.
\end{proof}

By Example \ref{ex:rickard}, the class of derived standardly stratified algebras is not closed under derived Morita equivalences.
It is therefore natural to ask which dg algebras are derived Morita equivalent to derived standardly stratified algebras.
By Theorem \ref{thm:dgqh}, every smooth proper derived standardly stratified algebra over an algebraically closed field is derived quasi-hereditary.
Combining this with Corollary \ref{cor:exc_char}, we obtain the following answer for smooth proper dg algebras.

\begin{corollary}\label{cor:smooth_exc}
Let $\kk$ be an algebraically closed field, and let $B$ be a smooth proper dg algebra.
Then the following statements are equivalent.
\begin{enumerate}
\item $B$ is derived Morita equivalent to a proper derived standardly stratified algebra.
\item $B$ is derived Morita equivalent to a proper derived quasi-hereditary algebra.
\item $\per B$ admits a full exceptional sequence.
\end{enumerate}
\end{corollary}
\begin{proof}
(2)$\Rightarrow$(1): This is clear.

(1)$\Rightarrow$(3): Let $A$ be a proper derived standardly stratified algebra which is derived Morita equivalent to $B$.
Since smoothness is invariant under derived Morita equivalences, $A$ is smooth.
Hence $\per A$ admits a full exceptional sequence by Theorem \ref{thm:dgqh}, and so does $\per B$.

(3)$\Rightarrow$(2): Let $\E_{B}$ be the dg category of homotopically projective dg $B$-modules belonging to $\per B$.
Since $B$ is proper, every full exceptional sequence in $\per B$ satisfies (HF).
Applying the proof of (2)$\Rightarrow$(1) in Corollary \ref{cor:exc_char} to $\E_{B}$, we obtain a proper derived quasi-hereditary algebra $A$ and a triangle equivalence $\per B\to\per A$ induced by a dg functor.
Hence $B$ is derived Morita equivalent to $A$.
\end{proof}

By Corollary \ref{cor:smooth_exc}, a smooth proper dg algebra $B$ such that $\per B$ has no exceptional objects is not derived Morita equivalent to any proper derived standardly stratified algebra.
The following example shows that such an algebra exists, even among finite-dimensional algebras.

\begin{example}\label{ex:happel}
Let $\kk$ be an algebraically closed field, and let $B=\kk Q/I$, where $Q$ is the quiver
\[
\begin{tikzcd}
1 \arrow[r, bend left=40, "a_{1}"] \arrow[r, "a_{2}"] & 2 \arrow[l, bend left=40, "b"]
\end{tikzcd}
\]
and $I=\langle a_{1}b,\,ba_{2}\rangle$.
The algebra $B$ has global dimension $3$, and hence it is smooth.
However, $\per B$ contains no exceptional objects (see \cite[Section 3]{H91} and \cite[Example 4.10]{RS22}).
By Corollary \ref{cor:smooth_exc}, $B$ is not derived Morita equivalent to any proper derived standardly stratified algebra, although it has finite global dimension.
\end{example}

Corollary \ref{cor:smooth_exc} does not apply to dg algebras which are not smooth, such as symmetric algebras which are not semisimple.
For such symmetric algebras, we have the following obstruction.

\begin{corollary}\label{cor:symmetric}
Let $B$ be a connected symmetric finite-dimensional algebra with at least two isomorphism classes of simple modules.
Then $\per B$ is not triangle equivalent to $\per A$ for any proper derived standardly stratified algebra $A$.
\end{corollary}
\begin{proof}
Suppose that $\per B$ is triangle equivalent to $\per A$ for some proper derived standardly stratified algebra $A$ with respect to a pair $(\PP,\leq)$ indexed by $\Lambda$.
Write $n=|\Lambda|$ and identify $\Lambda$ with $[1,n]$ via the order-preserving bijection.
By Remark \ref{rem:dss_indec} and \cite[Theorem 2.27]{AI12}, $n$ is the rank of $K_{0}(\per A)\cong K_{0}(\per B)$, that is, the number of isomorphism classes of simple $B$-modules.
In particular, $n\geq 2$.
By (DS2) and Lemma \ref{lem:thick=per}, we have a semiorthogonal decomposition 
\begin{align}
\per A=\langle\thick_{\D(A)}(\bde(1)),\ldots,\thick_{\D(A)}(\bde(n))\rangle,\notag
\end{align}
which gives a semiorthogonal decomposition $\per B=\langle\C_{1},\ldots,\C_{n}\rangle$ with nonzero components.
Since $B$ is symmetric, we have $\Hom_{\per B}(X,Y)\cong\Hom_{\per B}(Y,X)^{\ast}$ for all $X,Y\in\per B$ (see \cite{H88}).
Thus $\Hom_{\per B}(\C_{i},\C_{j})=0$ for all $i\neq j$, and hence $B\cong B_{1}\oplus\cdots\oplus B_{n}$ with $B_{i}\in\C_{i}$.
Since $B$ generates $\per B$, each $B_{i}$ is nonzero.
Therefore $B\cong\End_{\per B}(B)\cong\prod_{i=1}^{n}\End_{\per B}(B_{i})$ is a product of at least two nonzero algebras, which contradicts the connectedness of $B$.
\end{proof}

\subsection{A generalization to arbitrary dg algebras}

We mention a possible generalization of derived standardly stratified algebras to arbitrary dg algebras.
In this setting, we allow arbitrary shifts in condition (DS3).

\begin{definition}\label{def:gdss}
Let $(\Lambda,\leq)$ be a finite totally ordered set. Fix a dg algebra $A$ and a family $\PP=(P(i))_{i\in\Lambda}$ of objects in $\D(A)$ satisfying $\add\bigoplus_{i\in\Lambda}P(i)=\add A$.
The dg algebra $A$ is called a \emph{derived standardly stratified algebra} with respect to $(\PP,\leq)$ if there exist triangles in $\D(A)$
\begin{align}
K(i)\to P(i)\to \bde(i)\to K(i)[1]\quad(i\in\Lambda)\notag
\end{align}
satisfying (DS1) and (DS2) in Definition \ref{def:dss}, and the following condition:
\begin{enumerate}
\item[(DS3$'$)] $K(i)\in\F^{\ZZ}(\bde(>i)):=\F(\{\bde(j)[r]\mid j>i,\ r\in\ZZ\})$ for all $i\in\Lambda$.
\end{enumerate}
\end{definition}

When $A$ is connective, Definition \ref{def:gdss} is equivalent to Definition \ref{def:dss} for each fixed pair $(\PP,\leq)$ by the following proposition.

\begin{proposition}\label{prop:gdss_connective}
Let $A$ be a derived standardly stratified algebra with respect to $(\PP,\leq)$ in the sense of Definition \ref{def:gdss}.
If $A$ is connective, then $K(i)\in\F^{\geq 0}(\bde(>i))$ for all $i\in\Lambda$.
\end{proposition}
\begin{proof}
We show that $K(i)\in\F^{\geq 0}(\bde(>i))$ and $\bde(i)\in\D^{\leq 0}(A)$ by descending induction on $i$.
If $i$ is the maximal element of $\Lambda$, then $K(i)=0$ by (DS3$'$), and hence $\bde(i)\cong P(i)\in\D^{\leq 0}(A)$.

Assume that the assertions hold for all $j>i$.
Since $A$ is connective and $\bde(j)\in\D^{\leq 0}(A)$, we have $\Hom_{\D(A)}(P(j),\bde(j)[r])=0$ for all $j>i$ and $r>0$.
Moreover, (DS2) and (DS3$'$) imply $\Hom_{\D(A)}(K(j),\bde(j)[r])=0$ for all $j>i$ and $r\in\ZZ$.
Applying $\Hom_{\D(A)}(-,\bde(j))$ to the defining triangle $K(j)\to P(j)\to\bde(j)\to K(j)[1]$, we obtain $\Hom_{\D(A)}(\bde(j),\bde(j)[r])=0$ for all $j>i$ and $r>0$.
Thus $\bde(>i)$ is a higher standardizable set by (DS1) and (DS2).
In particular, $\F^{\geq 0}(\bde(>i))$ is closed under direct summands.

By descending induction using (DS3$'$), we have $\thick_{\D(A)}(Q_{i})=\thick_{\D(A)}(\bde(>i))$, where $Q_{i}=\bigoplus_{j>i}P(j)$.
Thus $K(i)\in\thick_{\D(A)}(Q_{i})$.
Moreover, (DS2) gives $\Hom_{\D(A)}(Q_{i},\bde(i)[r])=0$ for all $r\in\ZZ$.
Applying $\Hom_{\D(A)}(Q_{i},-)$ to the defining triangle for $i$, we obtain
\begin{align}
\Hom_{\D(A)}(Q_{i},K(i)[r])\cong\Hom_{\D(A)}(Q_{i},P(i)[r])=0\notag
\end{align}
for all $r>0$, where the last equality follows from the connectivity of $A$.
Since $A$ is connective and $Q_{i}\in\add A$, the object $Q_{i}$ is silting in $\thick_{\D(A)}(Q_{i})$.
Hence $K(i)\in\F^{\geq 0}(\add Q_{i})$ by \cite[Proposition 2.8(a) and Lemma 4.1]{IY18}.
On the other hand, for each $j>i$, the defining triangle for $j$ and the induction hypothesis give $P(j)\in\F^{\geq 0}(\bde(>i))$.
Since $\F^{\geq 0}(\bde(>i))$ is closed under direct summands, we obtain $\F^{\geq 0}(\add Q_{i})\subseteq\F^{\geq 0}(\bde(>i))$.
Therefore $K(i)\in\F^{\geq 0}(\bde(>i))$.

Finally, $K(i)\in\D^{\leq 0}(A)$ since $\bde(j)\in\D^{\leq 0}(A)$ for all $j>i$.
Hence the defining triangle for $i$ gives $\bde(i)\in\D^{\leq 0}(A)$.
This completes the induction.
\end{proof}

\subsection*{Acknowledgments}
TA is supported by JSPS Grant-in-Aid for Scientific Research JP24K06671.
MT is supported by JSPS Grant-in-Aid for Early-Career Scientists JP23K12959 and by the RDP Program Aiming at Maximizing the Career Potential of Female Researchers, Nagoya University (MEXT's Initiative for Realizing Diversity in the Research Environment, Leadership training type for women).

\subsection*{Use of artificial intelligence}
The authors used ChatGPT (OpenAI) and Claude (Anthropic) to assist in improving some of the mathematical arguments and the exposition, and in editing the English. 
The authors take full responsibility for the mathematical content and the final text.

\end{document}